\documentclass[pdflatex,sn-mathphys-num]{sn-jnl}

\usepackage{graphicx}%
\usepackage{amsmath,amssymb,amsfonts}%
\usepackage{amsthm}%
\usepackage[title]{appendix}%
\usepackage{xcolor}%
\usepackage{booktabs}%
\usepackage{enumitem}%

\theoremstyle{thmstyleone}%
\newtheorem{theorem}{Theorem}[section]%
\newtheorem{proposition}[theorem]{Proposition}%
\newtheorem{lemma}[theorem]{Lemma}%
\newtheorem{conjecture}[theorem]{Conjecture}%
\theoremstyle{thmstyletwo}%
\newtheorem{remark}[theorem]{Remark}%
\theoremstyle{thmstylethree}%

\newcommand{\R}{\mathbb{R}}
\newcommand{\PP}{\mathbb{P}}
\newcommand{\EE}{\mathbb{E}}
\newcommand{\LL}{\mathbb{L}}
\newcommand{\cE}{\mathcal{E}}
\newcommand{\cF}{\mathcal{F}}
\newcommand{\cD}{\mathcal{D}}
\newcommand{\cM}{\mathcal{M}}
\newcommand{\cL}{\mathcal{L}}
\newcommand{\cS}{\mathcal{S}}
\newcommand{\cZ}{\mathcal{Z}}
\newcommand{\cN}{\mathcal{N}}
\newcommand{\cG}{\mathcal{G}}
\newcommand{\cK}{\mathcal{K}}
\newcommand{\cO}{\mathcal{O}}
\newcommand{\cQ}{\mathcal{Q}}
\newcommand{\cU}{\mathcal{U}}
\newcommand{\cB}{\mathcal{B}}
\newcommand{\cH}{\mathcal{H}}
\newcommand{\kB}{k_{\mathrm{B}}}
\newcommand{\eps}{\varepsilon}
\newcommand{\Er}{E_{\mathrm{rest}}}
\newcommand{\Var}{\operatorname{Var}}
\newcommand{\Ent}{\operatorname{Ent}}
\newcommand{\capa}{\operatorname{Cap}}
\providecommand{\T}{}\renewcommand{\T}{\mathsf{T}}

\begin{document}

\title[A stochastic thermomechanical Cucker--Smale model]{On the stochastic thermomechanical
Cucker--Smale model with structure preserving noises}

\author[1]{\fnm{Seung-Yeal} \sur{Ha}}\email{syha@snu.ac.kr}
\author[2,3]{\fnm{Ruimeng} \sur{Hu}}\email{rhu@ucsb.edu}
\author[4,5]{\fnm{Tommaso} \sur{Ruggeri}}\email{tommaso.ruggeri@unibo.it}
\author*[2]{\fnm{Xu} \sur{Yang}}\email{xuyang@math.ucsb.edu}

\affil[1]{\orgdiv{Department of Mathematical Sciences and Research Institute of Mathematics}, \orgname{Seoul National University},
\orgaddress{\city{Seoul}, \country{Republic of Korea}}}
\affil[2]{\orgdiv{Department of Mathematics}, \orgname{University of California, Santa Barbara},
\orgaddress{\city{Santa Barbara}, \state{CA}, \country{USA}}}
\affil[3]{\orgdiv{Department of Statistics and Applied Probability}, \orgname{University of
California, Santa Barbara}, \orgaddress{\city{Santa Barbara}, \state{CA}, \country{USA}}}
\affil[4]{\orgdiv{Department of Mathematics}, \orgname{University of Bologna}, \orgaddress{\city{Bologna}, \country{Italy}}}
\affil[5]{\orgname{Accademia dei Lincei}, \orgaddress{\city{Rome}, \country{Italy}}}

\abstract{We construct a stochastic thermomechanical Cucker--Smale model from the entropy-based exchange structure of multivelocity, multitemperature gas mixtures. Noise acts on pairwise exchanges of momentum and total energy, preserving both conservation laws pathwise. For a prescribed Onsager mobility and reference measure, requiring reversibility of the Einstein distribution determines the noise covariance. With constant coefficients and connected exchange graphs, we establish a sharp necessary and sufficient condition for strictly positive constituent temperatures, expressed through the heat capacities of a coarse-graining cell. If it fails, at least one temperature reaches zero in finite time almost surely. In the positive regime, the Einstein distribution is the unique reversible equilibrium and relaxation is exponential, with a rate bound independent of the noise scale. The linearised stochastic and deterministic systems share the same relaxation rates, whose inverses are the classical momentum and thermal relaxation times for binary mixtures. We also prove Galilean covariance and identify the resulting correlation between momentum and energy fluctuations in moving frames.}


\maketitle

\section{Introduction}\label{sec:intro}

Collective behaviour appears in many biological systems, such as aggregation
of bacteria~\cite{BudreneBerg1991}, herding of sheep~\cite{Ginelli2015},
flocking of birds~\cite{Ballerini2008,CuckerSmale2007}, synchronised flashing of
fireflies and the activity of pacemaker cells and neurons~\cite{Buck1988,Winfree1967,MirolloStrogatz1990},
and schooling of fish~\cite{DegondMotsch2008}.
The development of mathematical models of collective motion was strongly
influenced by the work of Vicsek and collaborators~\cite{Vicsek1995},
followed by substantial work in statistical physics and
control~\cite{TonerTu1995,Jadbabaie2003,VicsekZafeiris2012} and by further models of
flocking~\cite{HaLiu2009,CarrilloFornasierRosadoToscani2010,MotschTadmor2011}.
The specific question addressed here is how this collective motion can be
connected to thermodynamic exchange processes.\par
\medskip

The classical Cucker--Smale model describes the emergence of collective velocity alignment through interactions among agents~\cite{CuckerSmale2007}. Its thermomechanical extension, introduced by Ha and Ruggeri~\cite{HaRuggeri2017}, has a specific physical origin: it was obtained by recognising the connection between the spatially homogeneous balance laws of multivelocity, multitemperature gas mixtures and the equations of collective motion. In this approach the temperature is not merely an additional variable assigned to an agent. The governing equations follow from the thermodynamic structure of mixture theory, including momentum and energy balances and the entropy principle; see also the treatment in~\cite[Chaps.~28 and 30]{RuggeriSugiyama2021}. The relation between the phenomenological description of gas mixtures, kinetic models and flocking was further investigated by Bae, Ha, Hwang and Ruggeri~\cite{BaeHaHwangRuggeri2024}.

The deterministic model couples mechanical alignment with thermal relaxation through exchanges among constituents. These exchanges conserve total momentum and total energy, while the thermodynamic entropy has nonnegative production. The resulting framework has motivated further developments, including a relativistic counterpart~\cite{HaKimRuggeri2020} and models of agents evolving on manifolds~\cite{HaKimSchloeder2021,AhnHaShim2021}. Outside the setting of a gas mixture, an agent may represent, for instance, a designed system with internal degrees of freedom. In such applications the physical meaning of the temperature-like variable must be specified rather than assumed.

An important question is how fluctuations can be introduced without losing these thermodynamic properties. Stochastic Cucker--Smale models with different noise mechanisms have been considered previously~\cite{HaLeeLevy2009,AhnHa2010,DuanFornasierToscani2010,BolleyCanizoCarrillo2011}, but perturbing the velocities of individual agents is not equivalent to allowing random internal exchanges. In particular, independently imposed velocity noise need not preserve the total momentum, total energy or strict positivity of the constituent temperatures. We therefore return to the mixture interpretation and make the exchange mechanisms, rather than the individual velocities, the point of departure.

This choice is also motivated by fluctuating thermodynamics. In classical fluctuating hydrodynamics, random stresses and heat fluxes describe fluctuations of transport processes~\cite{LandauLifshitzFM,LandauLifshitzSP,Balakrishnan2014}. By contrast, the production terms of rational extended thermodynamics describe local relaxation and remain meaningful even in spatially homogeneous mixtures; fluctuations associated with these terms have been studied, for example, in~\cite{Ikoma2011,RuggeriSugiyama2021}. Conservative stochastic exchanges are also familiar from energy-conserving dissipative particle dynamics~\cite{Espanol1997,BonetAvalosMackie1997}, and positivity issues in its numerical approximation have been considered in~\cite{Stoltz2017}. Our aim is not to claim that conservative noise is new, but to establish what its combination with the specific Ha--Ruggeri thermodynamic exchange structure implies for positivity, equilibrium and relaxation.

We consider a closed, spatially homogeneous system of an arbitrary number of constituents, with fixed partial densities, constant specific heats and constant symmetric momentum- and energy-exchange coefficients. The fluctuations occur in equal-and-opposite pairwise exchanges, so the total momentum and energy are conserved along each realisation. Once the dissipative mobility and the reference measure for the conserved variables have been chosen, requiring reversibility with respect to the Einstein distribution determines the constant noise covariance through a fluctuation--dissipation relation. This construction does not, by itself, provide a microscopic derivation of the noise or establish that the chosen reference measure is universal.

Our first main result gives a necessary and sufficient condition for all the temperatures to remain strictly positive. When the two exchange graphs are connected, positivity holds for every positive initial state precisely when the heat capacity of each constituent in the coarse-graining cell is at least the Boltzmann constant. The critical case is included, as are intersections of the boundaries where several temperatures vanish. If the condition is violated, the process reaches a zero-temperature boundary in finite time almost surely. Thus the positivity question has a sharp answer for the present constant-coefficient model, rather than merely a sufficient small-noise condition.

Under this positivity condition, the Einstein distribution is the unique invariant law on each manifold of fixed total momentum and energy. We establish exponential convergence in the appropriate square-integrable and relative-entropy classes, with lower bounds on the relaxation rate that are uniform in the fluctuation scale. The linearised stochastic dynamics has the deterministic relaxation rates; in the binary case, their inverses are the familiar momentum and thermal relaxation times of mixture theory~\cite[Sec.~28.6]{RuggeriSugiyama2021}. The corresponding limit of the spectral gap for the full nonlinear model remains open.

An additional issue, important in continuum thermodynamics, is Galilean covariance. The original deterministic equations can be written in an arbitrary inertial frame and then reduced to the barycentric rest frame. We show that the same reduction is legitimate for the stochastic model: conservation of total momentum makes the barycentric velocity constant along every trajectory. However, in a moving frame the random energy exchange must contain the work of the random momentum exchange. Consequently, the momentum and energy noises are correlated outside the rest frame, even though the rest-frame covariance has separate momentum and energy blocks.

Our analysis uses Dirichlet forms and the capacity of the zero-temperature boundary to handle positivity, and functional inequalities for relaxation. The present results concern a closed homogeneous model; spatial interactions, external forcing, physically derived temperature-dependent mobilities and a microscopic interpretation of the fluctuations remain separate questions.

The paper is organised as follows. Section~\ref{sec:model} presents the deterministic and stochastic models. Section~\ref{sec:main} states the main results. Sections~\ref{sec:pos} and~\ref{sec:relax} contain the proofs of positivity and relaxation, respectively. Galilean covariance is addressed in Section~\ref{sec:gal}, followed by a discussion of open questions in Section~\ref{sec:disc}. Technical details concerning Dirichlet forms are collected in Appendix~\ref{app:forms}.

\section{The model}\label{sec:model}
We first recall the deterministic thermomechanical exchanges and write them in Onsager form. Then we add internal noise and show how its covariance is determined.

\subsection{The deterministic model} \label{sec:2.1}
Consider a spatially homogeneous, non-reacting mixture of $n\ge2$ constituents in the Euclidean space ${\mathbb R}^d$. We assume that the constituent $a$ has a constant density $\rho_a>0$, a constant specific heat $c_{v,a}>0$, a velocity $v_a\in\R^d$ and a temperature $T_a>0$. We  also write 
\[ C_a=\rho_ac_{v,a}:~\mbox{heat capacity per unit volume} \quad \mbox{and} \quad \rho=\sum_a\rho_a:~\mbox{total(bulk) density}. \]
As state variables we use the momentum and energy densities of each constituent:
\begin{equation}\label{eq:state}
 p_a :=\rho_av_a\in\R^d,\quad E_a :=C_aT_a+\frac{|p_a|^2}{2\rho_a},\quad
 x:=(p_a,E_a)_{a=1}^n\in\R^{n(d+1)},
\end{equation}
where $|\cdot|$ is the standard $\ell^2$-norm in the Euclidean space ${\mathbb R}^d$. The temperature is then a function of the state:
\begin{equation} \label{eq:temp}
T_a(x)= \frac{E_a-|p_a|^2/(2\rho_a)}{C_a}. 
\end{equation}
We also use the following notations: for $a \in [n]$, 
\begin{equation}
\hspace{-0.2cm}
\begin{cases} \label{eq:bulk}
\displaystyle P:=\sum_ap_a: \mbox{ total momentum density},~~ E:=\sum_{a=1}^{n} E_a: \mbox{total energy density}, \\
\displaystyle v :=\frac{P}{\rho}: \mbox{barycentric velocity},~~ u_a :=v_a-v: \mbox{ diffusion velocities}.
\end{cases}
\end{equation}
Up to an additive constant, the entropy density of the mixture is given as follows.
\begin{equation}\label{eq:eta}
 \eta(x)=\sum_aC_a\log\frac{T_a(x)}{T_{\rm ref}},
\end{equation}
where $T_{\rm ref}>0$ is a fixed reference temperature. Then, it follows from \eqref{eq:state}, \eqref{eq:temp} and \eqref{eq:eta} that 
\[
 \frac{\partial\eta}{\partial E_a}=\frac1{T_a},\qquad
 \frac{\partial\eta}{\partial p_a}=-\frac{v_a}{T_a},
\]
The second derivatives are also given as follows. For variations $\xi_a\in\R^d$ of $p_a$ and $\theta_a\in\R$ of $E_a$,
\begin{equation}\label{eq:hessform}
 -\nabla^2\eta(x)[\xi,\theta]=\sum_a\Big[\frac{(\theta_a-v_a\cdot\xi_a)^2}{C_aT_a^2}
 +\frac{|\xi_a|^2}{\rho_aT_a}\Big].
\end{equation}
This follows from the following relations:
\[ C_a\nabla T_a=(-v_a,1) \quad \mbox{and} \quad C_a\nabla^2T_a=-\rho_a^{-1}(I_d\oplus0), \]
where the gradient and Hessian are taken in the variables $(p_a,E_a)$. In particular, the relation \eqref{eq:hessform} implies that $\eta$ is strictly concave. We assume that the constituents exchange momentum and energy with constant coefficients $\Lambda_{ab}=\Lambda_{ba}\ge0$ and $K_{ab}=K_{ba}\ge0$, $a\ne b$, and we take the following ansatz for the production terms based on entropy principle and main field theory (see ~\cite{HaRuggeri2017}):
\begin{equation}\label{eq:prod}
 F_a=\frac{1}{n} \sum_{b\ne a}\Lambda_{ab}\Big(\frac{u_b}{T_b}-\frac{u_a}{T_a}\Big),\qquad
 Q_a=\frac{1}{n} \sum_{b\ne a}K_{ab}\Big(\frac1{T_a}-\frac1{T_b}\Big), \quad a \in [n],
\end{equation}
so that 
\begin{equation} \label{eq:sum}
\sum_aF_a=0 \quad \mbox{and} \quad \sum_aQ_a=0. 
\end{equation}
In a general inertial frame, the deterministic model is given as follows:
\begin{equation}\label{eq:det}
\begin{cases}
\displaystyle \dot p_a=F_a, \quad a \in [n], \\
\displaystyle \dot E_a=Q_a+v\cdot F_a.
\end{cases}
\end{equation}
These are the momentum and total-energy balances in a generic inertial frame.  It follows from \eqref{eq:bulk}, \eqref{eq:sum} and \eqref{eq:det} that 
\[
\frac{dP}{dt} = \sum_{a} F_a = 0 \quad \mbox{and} \quad  \frac{dE}{dt} = \sum_{a} Q_a + v \cdot \sum_{a} F_a =  0,
\]
i.e., $P$ and $E$ are conserved along \eqref{eq:det}. Thus, the bulk velocity $v$ is constant. In the barycentric rest frame in which $v=0$, system reduces to 
\begin{equation} \label{eq:new-1}
\begin{cases}
\displaystyle \dot p_a=F_a, \quad a \in [n], \\
\displaystyle \dot E_a=Q_a
\end{cases}
\end{equation}
 with $u_a=v_a$. Moreover, with $\rho_a=c_{v,a}=1$ system \eqref{eq:new-1} is the thermomechanical Cucker--Smale model~\cite{HaRuggeri2017} with constant communication weights, where the coupling strengths $\kappa_1$ and $\kappa_2$ of~\cite{HaRuggeri2017} are absorbed into $\Lambda_{ab}$ and $K_{ab}$. \newline

We now write \eqref{eq:det} in Onsager form. Let $\cG_\Lambda$, $\cG_K$ be the weighted graph Laplacians:
\begin{equation} \label{eq:graph-laplacian}
(\cG_\Lambda)_{ab}= \begin{cases}
\displaystyle \sum_{c \ne a}\Lambda_{ac}, \quad & a = b, \\
\displaystyle -\Lambda_{ab}. \quad & a \neq b.
\end{cases}
\quad \mbox{and} \quad 
(\cG_K)_{ab}= \begin{cases}
\displaystyle \sum_{c \ne a} K_{ac}, \quad & a = b, \\
\displaystyle -K_{ab}. \quad & a \neq b,
\end{cases}
\end{equation}
and we define the matrix $ \LL$:
\begin{equation}\label{eq:L}
 \LL=L_{\mathrm m}\oplus L_{\mathrm e},\qquad
 L_{\mathrm m}=\frac1n\,\cG_\Lambda\otimes I_d,\qquad L_{\mathrm e}=\frac1n\,\cG_K ,
\end{equation}
acting on $x=\big((p_a)_a,(E_a)_a\big)$, and for $w\in\R^d$
\begin{equation}\label{eq:Lv}
 A_w=\begin{pmatrix}I_{nd}&0\\ B_w&I_n\end{pmatrix},\qquad (B_wz)_a=w\cdot z_a,\qquad
 \LL(w)=A_w\,\LL\,A_w^{\T} .
\end{equation}
Note that the matrix $\LL(w)$ is symmetric and positive semidefinite, and $\LL(0)=\LL$. By \eqref{eq:eta}, system \eqref{eq:det} becomes
\begin{equation}\label{eq:onsager}
 \dot x=\LL(v)\nabla\eta(x),\qquad
 \frac{d\eta}{dt}=\frac1n\sum_{a<b}\Big[\Lambda_{ab}\Big|\frac{u_a}{T_a}-\frac{u_b}{T_b}\Big|^2
 +K_{ab}\Big(\frac1{T_a}-\frac1{T_b}\Big)^2\Big]\ge0 .
\end{equation}
Note that $\LL$ is constant and does not depend on $\rho_a$ or $c_{v,a}$; these enter only through $\eta$. We prove the identity $\LL(v)\nabla\eta=(F_a,Q_a+v\cdot F_a)_a$ in Section~\ref{sec:gal}; in the rest frame it follows directly from \eqref{eq:eta}. The two exchange graphs are $G_\Lambda=\{(a, b):\Lambda_{ab}>0\}$ and $G_K=\{(a, b):K_{ab}>0\}$. 
\smallskip
\noindent\textbf{Condition (C).} Both graphs $G_\Lambda$ and $G_K$ are connected.
\smallskip

Under Condition~(C), the kernel of each graph Laplacian is spanned by constants.

\subsection{The stochastic model}
In this subsection, we present a stochastic version of the deterministic model \eqref{eq:onsager} by adding noise to the exchanges.  

Let $V>0$ and $\kB$ be the volume of a coarse-graining cell and the Boltzmann constant, respectively and we set
\begin{equation}\label{eq:eps}
 h=\frac{\kB}{V},\qquad \eps_a=\frac{h}{C_a}=\frac{\kB}{V\rho_ac_{v,a}} .
\end{equation}
For $a<b$ let $W_{ab}$ be independent standard Brownian motions in $\R^d$ and $\widetilde W_{ab}$ independent standard real Brownian motions, independent of the $W_{ab}$, and we set 
\[ W_{ba}=-W_{ab}, \quad \widetilde W_{ba}=-\widetilde W_{ab}. \]
The random exchanges are
\begin{equation}\label{eq:noise}
 dM_a=\sum_{b\ne a}\sqrt{\frac{2h\Lambda_{ab}}n}\;dW_{ab},\qquad
 dN_a=\sum_{b\ne a}\sqrt{\frac{2hK_{ab}}n}\;d\widetilde W_{ab},
\end{equation}
so that 
\[ \sum_adM_a=0 \quad \mbox{and} \quad \sum_adN_a=0. \]
Then, the proposed stochastic counterpart of the deterministic one \eqref{eq:det} is the following It\^o system:
\begin{equation}\label{eq:sde}
\begin{cases}
\displaystyle dp_a=F_a\,dt+dM_a, \quad a \in [n], \\
\displaystyle dE_a=\big(Q_a+v\cdot F_a\big)\,dt+dN_a+v\cdot dM_a. \\
\end{cases}
\end{equation}
Note that the term $v\cdot dM_a$ is the work done by the random momentum exchange in the barycentric motion. It is not a new source of noise. The noise in \eqref{eq:sde} has covariance $2h\LL(v)$, so the system has the same law as
\begin{equation}\label{eq:sdeL}
 dx=\LL(v)\nabla\eta\,dt+\sqrt{2h}\,A_v\LL^{1/2}\,d\beta, \quad  x =(p_a,E_a)_{a=1}^n,
\end{equation}
where $\beta$ is a standard Brownian motion in $\R^{n(d+1)}$. Total momentum and total energy are conserved along every path, so $v$ is constant. In the rest frame the model reads
\begin{equation}\label{eq:rest}
 dx=\LL\nabla\eta\,dt+\sqrt{2h\,\LL}\;d\beta .
\end{equation}
By Theorem~\ref{thm:gal} in the next section, every solution of \eqref{eq:sde} is a Galilean boost of a solution of \eqref{eq:rest}. For this reason we work in the rest frame until Section~\ref{sec:gal}. In the sequel, we discuss the state space and the Einstein density. \newline

\noindent $\bullet$ {\bf State space}:~For given $P \in\R^d$ and $E \in\R$ we set
\begin{align}
\begin{aligned} \label{eq:new-0}
 \cM_{P,E} &:=\Big\{x  \in {\mathbb R}^{n(d+1)}:\ \sum_ap_a(x) =P, \quad  \sum_aE_a=E \Big\},\\
 \cD_{P,E} &:=\{x\in\cM_{P,E}:\ T_a(x)>0,~~\forall~a \in [n] \},
\end{aligned}
\end{align}
and let $\Er=E-|P|^2/(2\rho)$ be the rest energy.

Note that the set $\cD_{P,E}$ is nonempty if and only if $\Er>0$. In the rest frame, we write
\[  \cM=\cM_{0,E} \quad \mbox{and} \quad \cD=\cD_{0,E}, \quad \mbox{with $E>0$}. \]
The set $\cD$ is a bounded convex open subset of the affine space $\cM$, which has dimension $m=(n-1)(d+1)$; it is convex because each $T_a$ is a concave function of $x$. Its boundary lies in the union of the faces $\partial_a\cD=\{x\in\overline\cD:T_a(x)=0\}$. We write $T\cM$ for the tangent space of $\cM$ and $dx$ for Lebesgue measure on $\cM$. Under Condition~(C), $\LL$ is positive definite on $T\cM$.

The coefficients of \eqref{eq:rest} are smooth in $\cD$, so for every initial state in $\cD$ there is a unique strong solution $X$ up to the exit time
\[
 \tau=\lim_{j\to\infty}\tau_j,\qquad
 \tau_j=\inf\{t\ge0:\ \operatorname{dist}(X_t,\partial\cD)\le 1/j\} .
\]
The temperatures stay positive for all time exactly when $\tau=\infty$. \newline

\noindent $\bullet$ {\bf The Einstein density:} ~On $\cD,$ we consider the probability density, with respect to $dx$,
\begin{equation}\label{eq:pi}
 \pi_h(x)=\frac1{\cZ_h}\,e^{\eta(x)/h}=\frac1{\cZ_h}\prod_a\Big(\frac{T_a(x)}{T_{\rm ref}}\Big)^{\alpha_a},
 \qquad \alpha_a=\frac1{\eps_a}=\frac{C_a}h .
\end{equation}
Since $V\eta$ is the entropy of a cell, this is the Einstein distribution $e^{V\eta/\kB}$ with respect to Lebesgue measure in the conserved variables, at fixed total momentum and energy.

The noise in \eqref{eq:rest} is not a free choice. Once we ask that the Einstein density has zero probability flux, its covariance is determined:
\begin{proposition}\label{prop:calib}
Let $\Sigma$ be a constant symmetric positive semidefinite matrix whose range lies in $T\cM$, and consider the diffusion process on $\cD$:
\[ dx=\LL\nabla\eta\,dt+\Sigma^{1/2}d\beta. \]
Then, the probability flux of the density $e^{\eta/h}$ vanishes on $\cD$ if and only if $\Sigma=2h\,\LL$.
\end{proposition}
\begin{proof}
The flux of a density $f$ is $J=\LL\nabla\eta\,f-\tfrac12\Sigma\nabla f$. Only the tangential gradient matters, because $\LL$ and $\Sigma$ are symmetric with range in $T\cM$, so they vanish on the normal space. For $f=e^{\eta/h}$ we get 
\[ J=\Delta\,\nabla\eta\,e^{\eta/h} \quad \mbox{with $\Delta=\LL-\Sigma/(2h)$}. \]
If $\Sigma=2h\LL$, then $J=0$. Conversely, we assume that 
\[ \Delta\nabla\eta(x)=0 \quad \mbox{for all $x\in\cD$}. \]
Differentiating along a tangent vector $\xi\in T\cM$ gives 
\[ \Delta\,\nabla^2\eta(x)\xi=0. \]
Let $\Pi$ be the orthogonal projection onto $T\cM$. Since $\Delta=\Delta\Pi$, the preceding identity gives 
\[ \Delta\,\Pi\nabla^2\eta(x)|_{T\cM}=0. \]
By \eqref{eq:hessform}, the operator $-\Pi\nabla^2\eta(x)|_{T\cM}$ is positive definite and therefore invertible. Thus $\Delta$ vanishes on $T\cM$ and on the normal space, so $\Delta=0$.
\end{proof}

This calibration has a limited meaning. Given $\LL$, the heat capacities, the cell volume $V$ and the reference measure, the fluctuation--dissipation relation fixes the covariance of the noise, and no extra noise parameter is left. It does not fix the factorisation of the covariance (the pair form \eqref{eq:noise} is one choice), and it does not derive the noise from a microscopic model. Note also that $VC_a$ is the heat capacity of constituent $a$ in a cell, so $\eps_a=\kB/(VC_a)$. For a classical gas whose molecules have mass $m_a$ and $D_a$ quadratic degrees of freedom, $c_{v,a}=D_a\kB/(2m_a)$, and with $\cN_a=\rho_aV/m_a$ molecules in a cell,
\begin{equation} \label{eq:epsgas}
 \eps_a=\frac{2}{D_a\cN_a} .
\end{equation}
Throughout the paper the coefficients $\rho_a$, $c_{v,a}$, $\Lambda_{ab}$ and $K_{ab}$ are constant, and the initial state is a fixed point of $\cD$ unless an initial distribution is given.

\section{Main results}\label{sec:main}
In this section, we list our results in the rest frame for the model \eqref{eq:rest} on $\cD$ with total energy $E>0$. Theorem~\ref{thm:gal} at the end of this section extends them to every inertial frame. The first result concerns the positivity of temperatures.

\begin{theorem}\label{thm:pos}
Under Condition~(C), the following assertions hold.
\begin{enumerate}[label=\textup{(\roman*)},leftmargin=2.2em]
\item If $\eps_a\le1$ for every $a$, then for every initial state in $\cD$  system \eqref{eq:rest} has a unique global strong solution, and almost surely $T_a(X_t)>0$ for all $t\ge0$ and all $a \in [n]$.
\item If $\eps_a>1$ for some $a$, then for every initial state in $\cD$, $\PP(\tau<\infty)=1$. The limit $X_{\tau-}$ exists almost surely and lies on $\partial\cD$, and almost surely $T_b(X_{\tau-})>0$ for every $b$ with $\eps_b\le1$.
\end{enumerate}
\end{theorem}
\begin{proof}
Since the proof need several preparatory estimates, we leave its proof in Section \ref{sec:4.3}.
\end{proof}
\begin{remark} \label{R3.2}
Note that under (C), all temperatures stay positive from every positive initial state if and only if $\max_a\eps_a\le1$. The condition depends only on the constituent with the smallest heat capacity in a cell; for a gas it reads $D_a\cN_a\ge2$. The assertion (ii) is about the first time the boundary is reached: at that time some constituent with $\eps_b>1$ has zero temperature, but not necessarily all of them.
\end{remark}

When the solution stays in $\cD$, we can ask its long-time behavior of \eqref{eq:rest}. Let $P_t$ be the transition semigroup of $X$, and for a probability density $f$ on $\cD$ let $H(f\mid\pi_h)=\int f\log(f/\pi_h)\,dx$ be the relative entropy. For $g\ge0$ we also write 
\[ \Ent_{\pi_h}(g)=\int g\log g\,\pi_h\,dx-\int g\,\pi_h\,dx\,\log\int g\,\pi_h\,dx, \]
Define the ensemble entropy of a probability density $f$ with $H(f\mid\pi_h)<\infty$ by
\begin{equation*}
 \cS_h[f]=\int_\cD f\,\eta\,dx-h\int_\cD f\log f\,dx=h\big(\log\cZ_h-H(f\mid\pi_h)\big) .
\end{equation*}

\begin{theorem}\label{thm:eq}
Assume that, for every initial state in $\cD$, the solution of \eqref{eq:rest} is global and remains in $\cD$ almost surely; by Theorem~\ref{thm:pos}(i), this holds under Condition~(C) when $\eps_a\le1$ for every $a$. Let $P_t$ denote its transition semigroup. Then the following assertions hold. 
\begin{enumerate}[label=\textup{(\roman*)},leftmargin=2.2em]
\item The Einstein density $\pi_h$ is invariant and reversible: for all bounded measurable functions $f$ and $g$, 
\[ \int g\,P_tf\,\pi_h\,dx=\int f\,P_tg\,\pi_h\,dx. \]
\item The measure $\pi_h\,dx$ is the unique invariant probability measure on $\cD$ if and only if the condition (C) holds. In that case $P_t(x,\cdot)\to\pi_h\,dx$ in total variation for every $x\in\cD$.
\item
If the initial law has a density $f_0$ with $H(f_0\mid\pi_h)<\infty$, and $f_t$ is the density of the law of $X_t$, then $t\mapsto\cS_h[f_t]$ is nondecreasing.
\end{enumerate}
\end{theorem}
\begin{proof}
We leave its proof in Section \ref{sec:5.2}.
\end{proof}
\begin{remark}
Uniqueness is meant on the state space $\cD$, where total momentum and energy are fixed. It is not enough that the union of the two graphs is connected: if one graph is disconnected, the momentum or the energy of each of its components is conserved separately. When $\eps_a>1$ for some $a$, we stop the solution at time $\tau$.
\end{remark}
Next, we study the speed of convergence using two quantities. The first one measures the curvature of the entropy relative to $\LL$. Let $\mathbb H(x)$ be the restriction of $-\nabla^2\eta(x)$ to $T\cM$, and we set
\begin{equation}\label{eq:kappa}
R :=\LL|_{T\cM}, \quad  \kappa :=\inf_{x\in\cD}\lambda_{\min}\big(R^{1/2}\mathbb H(x)R^{1/2}\big) .
\end{equation}
The second one comes from the linearization at the equilibrium. The entropy has its maximum on $\cD$ at the state $x^*$ with 
\[ v_a=0 \quad \mbox{and} \quad T_a=\bar T:=E/\sum_bC_b \quad \mbox{for all $a \in [n]$}. \]
With $D_\rho=\operatorname{diag}(\rho_a)$ and $D_C=\operatorname{diag}(C_a)$, let
\begin{equation}\label{eq:rates}
 \mu^\Lambda_k=\frac{\lambda_k\big(D_\rho^{-1/2}\cG_\Lambda D_\rho^{-1/2}\big)}{n\bar T},\qquad
 \mu^K_k=\frac{\lambda_k\big(D_C^{-1/2}\cG_KD_C^{-1/2}\big)}{n\bar T^2},\qquad k=2,\dots,n,
\end{equation}
where $0=\lambda_1\le\lambda_2\le\dots\le\lambda_n$ are the eigenvalues of the symmetric matrix in the bracket, and let $\lambda_*=\min(\mu^\Lambda_2,\mu^K_2)$. Finally, we need the Dirichlet form 
\[ \cE_h(f,f)=h\int_\cD\nabla f\cdot\LL\nabla f\,\pi_h\,dx \]
and its domain $\cF$, defined in Section~\ref{sec:frame}, and the spectral gap:
\begin{equation}\label{eq:gapdef}
 \lambda_{\rm gap}(h)=\inf\Big\{\frac{\cE_h(f,f)}{\Var_{\pi_h}(f)}:\ f\in\cF,\ \Var_{\pi_h}(f)>0\Big\}.
\end{equation}

\begin{theorem}\label{thm:relax}
Under Condition~(C), the following assertions hold.
\begin{enumerate}[label=\textup{(\roman*)},leftmargin=2.2em]
\item $\kappa>0$, and $\kappa$ depends only on $E$, $\rho_a$, $C_a$ and the exchange coefficients. In particular it does not depend on $h$.
\item For every $f\in\cF$,
\[
 \Var_{\pi_h}(f)\le\frac1\kappa\,\cE_h(f,f),\qquad
 \Ent_{\pi_h}(f^2)\le\frac2\kappa\,\cE_h(f,f) .
\]
\item If $\eps_a\le1$ for every $a$, then for every $f\in L^2(\pi_h)$ and every initial density $f_0$ with $H(f_0\mid\pi_h)<\infty$,  the following estimates hold.
\begin{gather*}
 \Big\|P_tf-\int f\,\pi_h\,dx\Big\|_{L^2(\pi_h)}\le e^{-\kappa t}\,
 \Big\|f-\int f\,\pi_h\,dx\Big\|_{L^2(\pi_h)},\\
 H(f_t\mid\pi_h)\le e^{-2\kappa t}\,H(f_0\mid\pi_h) .
\end{gather*}
\item The linearization of the deterministic system at $x^*$ has the relaxation rates \eqref{eq:rates}, where each $\mu^\Lambda_k$ has multiplicity $d$; the slowest rate is $\lambda_*>0$. The linearization of \eqref{eq:rest} at $x^*$ is an Ornstein--Uhlenbeck process, reversible with respect to the Gaussian law with covariance $h\mathbb H(x^*)^{-1}$, and its spectral gap is $\lambda_*$ for every $h>0$.
\item $\kappa\le\lambda_{\rm gap}(h)$ for all $h>0$, and $\limsup_{h\to0}\lambda_{\rm gap}(h)\le\lambda_*$. When $\eps_a\le1$ for every $a$, $\lambda_{\rm gap}(h)$ is the spectral gap of the semigroup $P_t$.
\end{enumerate}
\end{theorem}
\begin{proof} We leave the proof in Section \ref{sec:relax}.
\end{proof}
\begin{remark}
Assertions (i)--(iii) deal with the full nonlinear model. Assertion (iv) is only about the linearised model. For $n=2$ the rates in (iv) are
\[
 \mu^\Lambda_2=\frac{\Lambda_{12}}{2\bar T}\Big(\frac1{\rho_1}+\frac1{\rho_2}\Big),\qquad
 \mu^K_2=\frac{K_{12}}{2\bar T^2}\Big(\frac1{C_1}+\frac1{C_2}\Big).
\]
These are the inverse relaxation times $1/\tau_v$ and $1/\tau_T$ of the binary mixture in formula (28.33) of~\cite{RuggeriSugiyama2021}, with $\psi^0_{11}=\Lambda_{12}/2$ and $\theta^0_{11}=K_{12}/2$. The factor $1/2$ comes from the normalization $1/n$ in \eqref{eq:prod}, and for an ideal gas with ratio of specific heats $\Gamma_a=c_{p,a}/c_{v,a}$ one uses $c_{v,a}=\kB/(m_a(\Gamma_a-1))$. For identical constituents, the slowest rates are proportional to the algebraic connectivities of the two graphs. Assertions (iv) and (v) suggest the following conjecture, supported by the numerical example in Section~\ref{r:sec:num}.

\end{remark}
\begin{conjecture}\label{conj:gap}
Under the condition (C), the spectral gap converges to the relaxation rate $\lambda_*$ in the assertion (iv) in Theorem \ref{thm:relax}:
\[  \lim_{h \to 0} \lambda_{\rm gap}(h) = \lambda_*. \]
\end{conjecture}
The next theorem shows that restricting to the rest frame is no loss of generality.
\begin{theorem}\label{thm:gal}
For $c\in\R^d$ let $G_c$ be the Galilean boost
\[
 G_c:\ (p_a,E_a)_a\ \longmapsto\ \big(p_a+\rho_ac,\ E_a+c\cdot p_a+\tfrac12\rho_a|c|^2\big)_a .
\]
\begin{enumerate}[label=\textup{(\roman*)},leftmargin=2.2em]
\item 
$G_c$ is an affine bijection of $\cD_{P, E}$ onto $\cD_{P', E'}$, where $p'$ and $E^{\prime}$ are given by the following relations:
\[ P'=P+\rho c \quad \mbox{and} \quad E'=E+c\cdot P+\tfrac{1}{2} \rho |c|^2. \]
It preserves the temperatures $T_a$, the diffusion velocities $u_a$, the entropy $\eta$ and the rest energy $\Er$.
\vspace{0.1cm}
\item If $X$ solves \eqref{eq:sde}, then $G_c(X)$ solves \eqref{eq:sde} with the same Brownian motions. In particular, every solution of \eqref{eq:sde} is the image under $G_v$ of a solution of \eqref{eq:rest} with energy $\Er$.
\vspace{0.1cm}
\item The matrix field $x\mapsto\LL(v(x))$ is divergence free on $\R^{n(d+1)}$. Hence the probability flux of $e^{\eta/h}$ vanishes for \eqref{eq:sdeL} on the set $\{T_a>0\ \forall a\}$, and no extra drift is needed.
\vspace{0.1cm}
\item Theorem~\ref{thm:pos}, Theorem \ref{thm:eq} and~ Theorem \ref{thm:relax} hold on every $\cD_{P,E}$ with $\Er>0$, with $E$ replaced by $\Er$ and $\pi_h$ replaced by its image under $G_v$.
\end{enumerate}
\end{theorem}
\begin{proof} We leave the proof in Section \ref{sec:gal},
\end{proof}
\begin{remark}\label{rem:ref}
Below, we provide several comments on the results of Theorem \ref{thm:gal}. 
\begin{enumerate}
\item
In a frame with $v\ne0$ the momentum and energy noises are correlated, with covariance $2h\LL(v)$. The density $e^{\eta/h}$ cannot be normalised on the whole set $\{T_a>0\ \forall a\}$, because $\eta$ does not change under boosts; we normalise it at fixed total momentum and energy.
\vspace{0.1cm}
\item
Theorems~\ref{thm:pos}--\ref{thm:gal} remain true if $\cG_\Lambda/n$ and $\cG_K/n$ are replaced by any constant symmetric positive semidefinite $n\times n$ matrices whose kernels contain the constant vectors. Condition (C) is then replaced by the condition that both kernels consist only of the constant vectors, that is, that $\LL$ is positive definite on $T\cM$. This is the setting of Ruggeri and Sugiyama for non-reacting mixtures~\cite{RuggeriSugiyama2021}. The pair form \eqref{eq:noise} of the noise is available only when the off-diagonal entries are nonpositive.
\vspace{0.1cm}
\item
Suppose that the equilibrium law should be $\nu\,e^{\eta/h}$ with $\nu=\prod_aT_a^{\gamma_a}$. Keeping the same noise, the reversible model then has the extra drift $h\,\LL\nabla\log\nu$, so it is a different model. With $q_a=1/\eps_a+\gamma_a$, the equilibrium density can be normalised when $q_a>-1$ for all $a$. We expect that, when $q_a>0$, the arguments of Section~\ref{sec:pos} work with $\alpha_a$ replaced by $q_a$, so that the condition $\eps_a\le1$ in Theorem~\ref{thm:pos} becomes $q_a\ge1$. We do not pursue this.
\end{enumerate}
\end{remark}

\vspace{0.2cm}

In the following three sections, we provide proofs of Theorem \ref{thm:pos}, Theorem \ref{thm:eq}, Theorem \ref{thm:relax} and Theorem \ref{thm:gal} one by one. 

\section{Positivity of the temperatures}\label{sec:pos}
In this section we provide a proof of Theorem~\ref{thm:pos}. The proof strategy can be summarized as follows. We first introduce a Dirichlet form on $\overline\cD$ whose associated process agrees with $X$ before it reaches the boundary. A set of zero capacity is never reached by this process, while under (C) a set of positive capacity is reached almost surely. We show that the face has zero capacity when $\eps_a\le1$, and that part of it has positive capacity when $\eps_a>1$.

\subsection{Dirichlet forms and the process}\label{sec:frame}

For smooth functions $f$ and $g$, we set
\[
 \cE_h(f,g) :=h\int_\cD\nabla f\cdot\LL\nabla g\;\pi_h\,dx,\qquad
 \cE_{h,1}(f,f) :=\cE_h(f,f)+\|f\|^2_{L^2(\pi_h)},
\]
where $\nabla$ is the gradient along $\cM$. We use two closed extensions of this form,
\begin{equation}\label{eq:domains}
 \cF=\overline{C^\infty(\overline\cD)}^{\,\cE_{h,1}},\qquad
 \cF_0=\overline{C_c^\infty(\cD)}^{\,\cE_{h,1}}\subset\cF ,
\end{equation}
where $C^\infty(\overline\cD)$ denotes the restrictions to $\overline\cD$ of smooth functions on $\cM$. The generator of the model is
\begin{equation}\label{f:eq:gen}
 \cL f=\LL\nabla\eta\cdot\nabla f+h\operatorname{tr}(\LL\nabla^2f)=\frac{h}{\pi_h}\nabla\cdot\big(\pi_h\LL\nabla f\big) ;
\end{equation}
the second form uses $h\nabla\pi_h=\pi_h\nabla\eta$ and the fact that $\LL$ is constant. We collect the properties we need in three lemmas. Their proofs are standard and are given in Appendix~\ref{app:forms}. In the proofs, $\lambda_{\max}$ is the largest eigenvalue of $\LL$, and under (C), $\lambda^+_{\min}$ is its smallest eigenvalue on $T\cM$; norms are taken in the coordinates $(p_a,E_a)$. In the sequel, we state three preparatory lemmas whose proofs will be given in Appendix.

\begin{lemma}\label{f:lem:form}
\begin{enumerate}[label=\textup{(\roman*)},leftmargin=2.2em]
\item The form $\cE_h$ is closable on $C^\infty(\overline\cD)$ and on $C_c^\infty(\cD)$. The closures $(\cE_h,\cF)$ on $L^2(\overline\cD,\pi_h)$ and $(\cE_h,\cF_0)$ on $L^2(\cD,\pi_h)$ are regular, strongly local Dirichlet forms. The form $(\cE_h,\cF)$ is recurrent, hence conservative.
\vspace{0.1cm}
\item Every function that is Lipschitz on $\overline\cD$ belongs to $\cF$, and $\cE_h$ is given by the formula above with the almost everywhere gradient. If the function also has compact support in $\cD$, it belongs to $\cF_0$.
\vspace{0.1cm}
\item For $f\in C_c^\infty(\cD)$ and $g\in\cF$, $\cE_h(f,g)=-\int_\cD(\cL f)\,g\,\pi_h\,dx$.
\vspace{0.1cm}
\item Under (C), every $f\in\cF$ has a weak gradient in $L^2_{\rm loc}(\cD)$, and $\cE_h(f,g)=h\int_\cD\nabla f\cdot\LL\nabla g\,\pi_h\,dx$ for $f,g\in\cF$.
\end{enumerate}
\end{lemma}

By Lemma~\ref{f:lem:form} and~\cite{FOT}, there is a $\pi_h$-symmetric diffusion $Y$ on $\overline\cD$ associated with $(\cE_h,\cF)$, with infinite lifetime, and a $\pi_h$-symmetric diffusion $Y^0$ on $\cD$ associated with $(\cE_h,\cF_0)$, with lifetime $\zeta\le\infty$. Both have continuous paths and are defined up to an exceptional set of starting points. The process $Y$ provides a conservative extension on $\overline\cD$, while $Y^0$ is killed at the boundary. Let $\sigma_j$ and $\sigma^0_j$ be defined like $\tau_j$, with $Y$ and $Y^0$ in place of $X$, and let $\sigma$ be the hitting time of $\partial\cD$ by $Y$. The next lemma says that $X$, $Y$ and $Y^0$ agree until they come close to the boundary.

\begin{lemma}\label{f:lem:ident}
For quasi-every $y\in\cD$, the families of stopped processes $(X_{\cdot\wedge\tau_j})_{j\ge1}$, $(Y_{\cdot\wedge\sigma_j})_{j\ge1}$ and $(Y^0_{\cdot\wedge\sigma^0_j})_{j\ge1}$ have the same joint law under $\PP_y$. Consequently $\PP_y(\tau\le t)=\PP_y(\sigma\le t)=\PP_y(\zeta\le t)$ for all $t\ge0$.
\end{lemma}

The last lemma lets us pass from quasi-every starting point to every starting point.

\begin{lemma}\label{f:lem:ac}
Suppose Condition (C) holds. For every $x\in\cD$ and $s>0$, the law of $X_s$ restricted to $\{s<\tau\}$ is absolutely continuous with respect to $dx$. Moreover, for every smooth connected domain $\cD'\Subset\cD$, the solution killed on leaving $\cD'$ has a transition density that is continuous and strictly positive on $\cD'\times\cD'$.
\end{lemma}

\subsection{Capacity of the faces}
We now estimate the capacity of the faces. Capacities refer to the form $(\cE_h,\cF)$. For a relatively open set $\cU\subset\overline\cD$, $\capa(\cU)=\inf\{\cE_{h,1}(u,u):\ u\in\cF,\ u\ge1\ \pi_h\text{-a.e.\ on }\cU\}$, and this is extended to arbitrary sets by outer regularity. Since $\pi_h(N)\le\capa(N)$, sets of zero capacity are Lebesgue-null in $\cD$. In the next lemma, we use the cut-off functions $\varphi_r$, $0<r<1$:
\[
\varphi_r(t)= \begin{cases}
1, \quad & t\le r^2, \\
\frac{\log(r/t)}{\log(1/r)}, \quad & r^2\le t\le r, \\
0, \quad & t\ge r.
\end{cases}
\]
In what follows, temperatures are measured in units of $T_{\rm ref}$. We also use \eqref{eq:temp} and \eqref{eq:new-0} to see the bounds
\begin{equation}\label{f:eq:Tmax}
 T_a<T_a^{\max}:=E/C_a,\qquad |v_a|^2<2E/\rho_a\qquad\text{on }\cD,
\end{equation}
which hold because $0<C_aT_a\le E_a<E$ and $|p_a|^2/(2\rho_a)<E$. Next, we define
\[
 f_{a,r} :=\varphi_r(T_a) 
\]
and study its properties in the following lemma.
\begin{lemma}\label{f:lem:cap}
Suppose $\eps_a\le1$. Then, the following assertions hold. \newline
\begin{enumerate}[label=\textup{(\roman*)},leftmargin=2.2em]
\item $f_{a,r}$ is Lipschitz on $\overline\cD$, equals $1$ on the relatively open set $\{T_a<r^2\}\supset\partial_a\cD$.
\vspace{0.1cm}
\item
There exist positive constants $\Theta_a$ and $c_\pi$ such that 
\begin{equation}\label{f:eq:capbound}
 \cE_h(f_{a,r},f_{a,r})\le\frac{\Theta_a}{\log(1/r)},\qquad
 \int_\cD f_{a,r}^2\,\pi_h\,dx\le c_\pi|\cD|\,r^{\alpha_a}.
\end{equation}
\item
The face $\partial_a\cD$ has zero capacity:
\[ \capa(\partial_a\cD)=0. \]
\end{enumerate}
\end{lemma}
\begin{proof}
(i) First we collect three bounds on $\cD$. By \eqref{f:eq:Tmax}, $\pi_h\le c_\pi\,T_a^{\alpha_a}$ with $c_\pi:=\cZ_h^{-1}\prod_{b\ne a}(T_b^{\max})^{\alpha_b}$. The gradient of $T_a$ in the coordinates $(p_a,E_a)$ is $C_a^{-1}(-v_a,1)$, and its projection onto $\cM$ is not longer, so we have
\[ |\nabla T_a|\le g_a:=C_a^{-1}(1+2E/\rho_a)^{1/2}. \]
Finally, since $T_a$ is concave, $\{T_a>t\}\cap\cD$ is a convex subset of $\cD$, and the level set $\{T_a=t\}\cap\cD$ lies in its boundary. Surface area increases with convex sets, so we have
\[ \cH^{m-1}(\{T_a=t\}\cap\cD)\le P_\cD:=\cH^{m-1}(\partial\cD)<\infty. \]
The function $f_{a,r}$ is Lipschitz on $\overline\cD$, so it is in $\cF$ by Lemma~\ref{f:lem:form}(ii). \newline

\noindent (ii)~$\bullet$~(Derivation of the first estimate in \eqref{f:eq:capbound}):~We use 
\[ \nabla f\cdot\LL\nabla f\le\lambda_{\max}\varphi_r'(T_a)^2|\nabla T_a|^2 \le\lambda_{\max}g_a\varphi_r'(T_a)^2|\nabla T_a| \]
and the co-area formula on $\cM$ to get 
\begin{equation} \label{eq:new-2}
\hspace{-0.2cm} \cE_h(f_{a,r},f_{a,r})\le h\lambda_{\max}g_a\int_0^\infty\varphi_r'(t)^2
 \Big(\int_{\{T_a=t\}\cap\cD}\pi_h\,d\cH^{m-1}\Big)dt
 \le\Theta_a\int_{r^2}^r\varphi_r'(t)^2\,t^{\alpha_a}\,dt,
\end{equation}
with $\Theta_a:=h\lambda_{\max}g_ac_\pi P_\cD$.  \newline

\noindent On the interval $(r^2,r)$ we use  $t<1$ and $\alpha_a\ge1$ to see
\begin{equation} \label{eq:new-3}
\varphi_r'(t)=-1/(t\log(1/r)) \quad \mbox{and} \quad t^{\alpha_a-1}\le1.
\end{equation}
Now, we combine \eqref{eq:new-2} and \eqref{eq:new-3} to get
\[
 \int_{r^2}^r\varphi_r'(t)^2\,t^{\alpha_a}\,dt\le\frac1{\log^2(1/r)}\int_{r^2}^r\frac{dt}t
 =\frac1{\log(1/r)} .
\]
\noindent  $\bullet$~(Derivation of the second estimate in \eqref{f:eq:capbound}):~The second bound in \eqref{f:eq:capbound} holds because $f_{a,r}=0$ where $T_a\ge r$. \newline

\noindent (iii)~Since $f_{a,r}=1$ on $\{T_a<r^2\}$, we get 
\[ \capa(\partial_a\cD)\le\Theta_a/\log(1/r)+c_\pi|\cD|r^{\alpha_a}. \]
Since the right-hand side of the above relation tends to zero as $r \to 0$, we have the desired estimate.  Note that this is the only place where $\eps_a\le1$ is used, including the case $\eps_a=1$.
\end{proof}
\begin{remark} Below, we comment on the results of Lemma~\ref{f:lem:cap}.
\begin{enumerate}
\item
In one dimension, Lemma~\ref{f:lem:cap} is the classical fact that the point $0$ is reached by the diffusion with speed measure $t^{\alpha_a}dt$ and scale density $t^{-\alpha_a}$ exactly when $\int_0t^{-\alpha_a}dt<\infty$, that is, when $\alpha_a<1$.
\vspace{0.1cm}
\item
The cut-off functions of Lemma~\ref{f:lem:cap} also show that the two forms coincide when every face has zero capacity. We use this in Section~\ref{sec:relax}.
\end{enumerate}
\end{remark}

\begin{lemma}\label{f:lem:cap-1}
Suppose the condition (C) holds. If $\eps_a>1$, then $\partial_a\cD$ contains a compact set of positive capacity.
\end{lemma}
\begin{proof}
\noindent For some $a \in [n]$, we assume that $\eps_a > 1$. Then, we choose $x_0\in\partial_a\cD$ with $T_b(x_0)>0$ for $b\ne a$, for example with all momenta zero, $E_a=0$ and the energy shared among the other constituents. The tangential gradient of $T_a$ at $x_0$ is not zero, because moving energy from another constituent to $a$ changes $T_a$. So there is a smooth chart $x\mapsto(t,w)=(T_a(x),w(x))$ near $x_0$, with $w$ in a relatively compact open set $\cO\Subset\R^{m-1}$ and $t_0\in(0,1)$, such that the set $\cQ$ corresponding to $(0,t_0)\times\cO$ lies in $\cD$, the set $\cB$ corresponding to $\{0\}\times\overline\cO$ is a compact subset of $\partial_a\cD$, and the other temperatures are bounded below on $\cQ$. On $\cQ$, we then have
\[
 \pi_h\ge c_-t^{\alpha_a},\qquad dx\ge j_-\,dt\,dw,\qquad |\partial_tu|\le g_+|\nabla u|,
\]
for some constants $c_-,j_-,g_+>0$, and $\nabla u\cdot\LL\nabla u\ge\lambda^+_{\min}|\nabla u|^2$ by (C). Hence, for $u\in C^\infty(\overline\cD)$,
\begin{equation}\label{f:eq:loc}
 \int_\cO\!\int_0^{t_0}|\partial_tu|^2t^{\alpha_a}\,dt\,dw\le b_1\,\cE_h(u,u),\qquad
 \int_\cO\!\int_{t_0/2}^{t_0}u^2\,dt\,dw\le b_2\,\|u\|^2_{L^2(\pi_h)} .
\end{equation}
Since $\alpha_a<1$, the integral $I_a:=\int_0^{t_0}r^{-\alpha_a}dr$ is finite. For $w\in\cO$ and $0<s<t\le t_0$, the Cauchy--Schwarz inequality yields
\[
 |u(s,w)|\le|u(t,w)|+I_a^{1/2}\Big(\int_0^{t_0}|\partial_ru(r,w)|^2r^{\alpha_a}dr\Big)^{1/2}.
\]
Fix $\delta<t_0/2$, we take an average of the above relation over $s\in(0,\delta)$ and $t\in(t_0/2,t_0)$, integrate the resulting relation over $\cO$, and use the Cauchy--Schwarz inequality and \eqref{f:eq:loc} to find
\begin{equation}\label{f:eq:capbd}
 \frac1\delta\int_\cO\!\int_0^\delta|u|\,ds\,dw\le b_3\,\cE_{h,1}(u,u)^{1/2},
\end{equation}
with $b_3$ independent of $\delta$ and $u$. The left-hand side is continuous in the $\cE_{h,1}$-norm, because 
\[ \int_0^\delta|u|\,ds\le(\int_0^\delta u^2s^{\alpha_a}ds)^{1/2}(\int_0^\delta s^{-\alpha_a}ds)^{1/2}, \]
so \eqref{f:eq:capbd} holds for all $u\in\cF$. Now let $u\in\cF$ with $u\ge1$ a.e.\ on a relatively open set $\cU\supset\cB$. By compactness, $\cU$ contains the set corresponding to $(0,\delta)\times\cO$ for some $\delta>0$, so the left-hand side of \eqref{f:eq:capbd} is at least $|\cO|$. Thus we have
\[ \capa(\cB)\ge|\cO|^2/b_3^2>0. \]
\end{proof}
\begin{lemma}\label{f:lem:F0F}
If $\eps_a\le1$ for every $a \in [n]$, then we have 
\[ \cF_0=\cF. \]
\end{lemma}

\begin{proof}
Since $\cF_0$ is closed and $C^\infty(\overline\cD)$ is dense in $\cF$, it suffices to show that $u\in\cF_0$ for $u\in C^\infty(\overline\cD)$. We set 
\[ \chi_r :=\max_{a \in [n]} f_{a,r}. \]
Then, it is Lipschitz continuous on $\overline\cD$, and equals $1$ where some $T_a<r^2$, and $1-\chi_r$ vanishes outside the compact set $\{x\in\cD:T_a(x)\ge r^2\ \forall a\}\subset\cD$. Again, we write
\[ u=u(1-\chi_r)+u\chi_r. \]
 The first term is Lipschitz with compact support in $\cD$, so it is in $\cF_0$ by Lemma~\ref{f:lem:form}(ii). For the second term, a.e.\ $\nabla\chi_r=\nabla f_{a,r}$ for some $a$, so we have 
 \begin{align*}
 \begin{aligned}
&  \nabla\chi_r\cdot\LL\nabla\chi_r\le\sum_a\nabla f_{a,r}\cdot\LL\nabla f_{a,r} \quad \mbox{and} \\
& \cE_h(u\chi_r,u\chi_r)\le2h\lambda_{\max}\|\nabla u\|_\infty^2\int_\cD\chi_r^2\,\pi_h\,dx
 +2\|u\|_\infty^2\sum_a\cE_h(f_{a,r},f_{a,r}) .
\end{aligned}
\end{align*}
By \eqref{f:eq:capbound} and $\chi_r^2\le\sum_af_{a,r}^2$, both terms and $\|u\chi_r\|_{L^2(\pi_h)}$ tend to zero as $r\to0$. Hence $u(1-\chi_r)\to u$ in $\cE_{h,1}$, and $u\in\cF_0$.
\end{proof}

\subsection{Proof of Theorem~\ref{thm:pos}}  \label{sec:4.3}

Now, we are ready to provide the proof of our first main result on the positivity of temperature. In the sequel, we show each assertion one by one. \newline

\noindent $\bullet$~(Proof of the first assertion (i)):~Suppose that the condition (C) and  $\eps_a\le1$ for every $a \in [n]$ hold. \newline

$\diamond$~Step A.1 (quasi-every starting point):~By Lemma~\ref{f:lem:cap}(iii),  we have
\[ \capa(\partial\cD)\le\sum_a\capa(\partial_a\cD)=0. \]
For a regular Dirichlet form, sets of zero capacity are exceptional~\cite{FOT}, so $Y$ started from quasi-every $y\in\cD$ never reaches $\partial\cD$. By Lemma~\ref{f:lem:ident}, the function $\omega(y):=\PP_y(\tau<\infty)$ vanishes for quasi-every $y$, hence almost everywhere.

\vspace{0.1cm}

$\diamond$~Step A.2 (every starting point):~Fix $x\in\cD$ and $s>0$. By the Markov property and Lemma~\ref{f:lem:ac},
\[
 \PP_x(s<\tau<\infty)=\EE_x\big[\mathbf 1_{\{s<\tau\}}\,\omega(X_s)\big]=0 .
\]
Since $x$ is an interior point and the paths are continuous, $\PP_x(\tau\le s)\to0$ as $s\downarrow0$, so $\PP_x(\tau<\infty)=0$. Pathwise uniqueness up to $\tau$ gives a unique global strong solution, and all temperatures stay positive. 

Note that the argument does not distinguish faces from corners: a corner, where several temperatures vanish, lies in the union of the faces, which has zero capacity. \newline

\noindent $\bullet$~(Proof of the second assertion (ii)):~Assume (C) and $\eps_a>1$ for some $a$. \newline

$\diamond$~Step B.1 (irreducibility):~The form $(\cE_h,\cF)$ is recurrent by Lemma~\ref{f:lem:form}(i). If $\mathcal I$ is an invariant set, then we have
\[ 1_{\mathcal I}\in\cF \quad \mbox{and} \quad \cE_h(1_{\mathcal I},1_{\mathcal I})\le\cE_h(1,1)=0, \]
see~\cite{FOT}. By Lemma~\ref{f:lem:form}(iv), $1_{\mathcal I}$ has weak gradient zero in the connected set $\cD$, so it is almost everywhere constant. Hence $(\cE_h,\cF)$ is irreducible.

\vspace{0.1cm}

$\diamond$~Step B.2 (quasi-every starting point):~By Lemma~\ref{f:lem:cap-1} there is a compact set $\cB\subset\partial_a\cD$ of positive capacity. For an irreducible recurrent Dirichlet form, a Borel set of positive capacity is hit almost surely from quasi-every starting point~\cite[Theorem~4.7.1(iii)]{FOT}. So $\PP_y(\sigma<\infty)=1$ for quasi-every $y$, and by Lemma~\ref{f:lem:ident}, $\omega(y)=1$ for quasi-every $y$, so almost everywhere.

\vspace{0.1cm}

$\diamond$~Step B.3 (every starting point):~Fix $x\in\cD$ and $s>0$. As before,
\[
 \PP_x(\tau<\infty)\ge\EE_x\big[\mathbf 1_{\{s<\tau\}}\,\omega(X_s)\big]=\PP_x(s<\tau)
 \longrightarrow1\qquad(s\downarrow0).
\]

$\diamond$~Step B.4 (the boundary limit):~By Lemma~\ref{f:lem:ident}, for starting points outside an exceptional set, the families $(\tau_j,X_{\cdot\wedge\tau_j})_j$ and $(\sigma_j,Y_{\cdot\wedge\sigma_j})_j$ have the same joint law. On $\{\sigma<\infty\}$, which has full probability, $\sigma_j\uparrow\sigma$ and $Y_{\cdot\wedge\sigma_j}\to Y_{\cdot\wedge\sigma}$ locally uniformly in time, since $Y$ has continuous paths. So the same holds for $X$: the limit $X_{\tau-}$ exists, lies on $\partial\cD$, and has the law of $Y_\sigma$. \newline

Let $\mathcal V=\bigcup_{b:\,\eps_b\le1}\partial_b\cD$, which has zero capacity by Lemma~\ref{f:lem:cap}(iii), so that $Y_\sigma\notin\mathcal V$ almost surely for quasi-every start. Let $\mathcal A$ be the event that $\tau=\infty$, or that $\tau<\infty$ and $X_{\tau-}$ does not exist or does not lie in $\partial\cD\setminus\mathcal V$. Then we have
\[ \omega_{\mathcal A}(y):=\PP_y(\mathcal A)=0 \quad \mbox{for quasi-every $y$, hence almost everywhere}. \]
On $\{s<\tau\}$ the event $\mathcal A$ depends only on the path after time $s$, so the Markov property and Lemma~\ref{f:lem:ac} give 
\[ \PP_x(\mathcal A,\ s<\tau)=\EE_x[\mathbf 1_{\{s<\tau\}}\omega_{\mathcal A}(X_s)]=0. \]
Letting $s\downarrow0$, we have
\[ \PP_x(\mathcal A)=0 \quad \mbox{for every $x\in\cD$}. \]
\begin{remark}\label{f:rem:bessel}
The threshold $\eps_a=1$ can be understood as follows. Near a point of the face $\{T_a=0\}$ where the other temperatures are positive, $T_a$ has drift $(K_a+\Lambda_a|v_a|^2)/(nC_aT_a)+O(1)$ and quadratic variation rate $2h(K_a+\Lambda_a|v_a|^2)/(nC_a^2)$, with $\Lambda_a=\sum_{b\ne a}\Lambda_{ab}$ and $K_a=\sum_{b\ne a}K_{ab}$. So $T_a$ behaves locally like a Bessel process of dimension $1+1/\eps_a$, which avoids zero if and only if this dimension is at least $2$, that is, $\eps_a\le1$. The factor $K_a+\Lambda_a|v_a|^2$ is the same in the drift and in the noise because of the fluctuation--dissipation relation. This heuristic is not used in the proofs.
\end{remark}

\section{Equilibrium and relaxation}\label{sec:relax}

In this section, we provide proofs for Theorem~\ref{thm:eq} and Theorem \ref{thm:relax}. First, we identify the semigroup $P_t$ with the semigroup of the Dirichlet form and prove Theorem~\ref{thm:eq}, and then we prove functional inequalities for $\pi_h$, which yields the exponential decay. Finally we study the linearized model and the spectral gap for small noise.

\subsection{Equilibrium} \label{sec:5.1}

In this subsection, we assume, as in Theorem~\ref{thm:eq}, that for every initial state in $\cD$ the solution stays in $\cD$ for all time. Let $(P^{\cF}_t)$ and $(P^0_t)$ be the $L^2(\pi_h)$-semigroups of $(\cE_h,\cF)$ and $(\cE_h,\cF_0)$, respectively. 

\begin{proposition}\label{f:prop:semigroup}
\begin{enumerate}[label=\textup{(\roman*)},leftmargin=2.2em]
\item The process $Y^0$ is conservative, its law coincides with that of $X$ for quasi-every start, and $P^0_tf=P_tf$ $\pi_h$-almost everywhere for bounded measurable $f$.
\vspace{0.1cm}
\item If in addition $\eps_a\le1$ for every $a$, then $\cF_0=\cF$ and $P^{\cF}_t=P^0_t=P_t$.
\vspace{0.1cm}
\item If the initial law has a density $f_0$, then the law at time $t$ has the density $f_t=\pi_h\,P_t(f_0/\pi_h)$.
\end{enumerate}
\end{proposition}
\begin{proof}
\noindent (i) By Lemma~\ref{f:lem:ident}, we have
\[ \PP_y(\zeta\le t)=\PP_y(\tau\le t)=0 \quad \mbox{for quasi-every $y$}, \]
so $Y^0$ is conservative. Its stopped laws agree with those of $X$ at every level $j$, and $\tau_j\uparrow\infty$, so $Y^0$ and $X$ have the same law for quasi-every start. Exceptional sets are $\pi_h$-null, so we have
\[ P^0_tf=P_tf \quad \mbox{$\pi_h$-a.e.} \]

\noindent (ii) This follows from Lemma~\ref{f:lem:F0F}.

\vspace{0.1cm}

\noindent (iii) By (i), $P_t$ is symmetric in $L^2(\pi_h)$ on bounded functions, and by truncation and monotone convergence 
\[ \int g\,P_tf\,\pi_h\,dx=\int f\,P_tg\,\pi_h\,dx \quad \mbox{for all measurable $f,g\ge0$}. \]
For bounded measurable $\varphi$, this implies
\[ \int\varphi f_t\,dx=\int P_t\varphi\,f_0\,dx=\int\varphi\,P_t(f_0/\pi_h)\,\pi_h\,dx. \]
\end{proof}

\subsection{Proof of Theorem \ref{thm:eq}} \label{sec:5.2}
In this subsection, we provide a detailed proof for each assertion on the properties of the Einstein density  $\pi_h$ one by one. 

\vspace{0.2cm}

\noindent (i)~(First assertion):~By Proposition~\ref{f:prop:semigroup}(i), we have
\[ \int g\,P_tf\,\pi_h\,dx=\int f\,P_tg\,\pi_h\,dx \]
for bounded measurable $f,g$. This is reversibility, and taking $g=1$ gives invariance.

\vspace{0.2cm}

\noindent (ii)~(Second assertion):~We split the proof into two steps. \newline

\noindent $\bullet$ Step A ($\Longleftarrow$):~Suppose that (C) holds. Since $\cD$ is convex, there are smooth connected domains $\cD_k\Subset\cD$ increasing to $\cD$. By Lemma~\ref{f:lem:ac}, for $x\in\cD_k$ and $t>0$ we have 
\[ P_t(x,\cU)\ge\PP_x(X_t\in\cU,\ t<\tau_{\cD_k})>0 \]
for every nonempty open set $\cU\subset\cD_k$, where $\tau_{\cD_k}$ is the exit time from $\cD_k$. So $(P_t)$ is irreducible. The semigroup is also strong Feller. Let $f$ be bounded and measurable, let $\cK\subset\cD_k$ be compact, $\delta=\operatorname{dist}(\cK,\partial\cD_k)$, and let $b$ be a bound for the drift on $\cD_k$. For $x\in\cK$ and $0<s<t$ with $bs<\delta/2$, the Markov property gives
\[
 P_tf(x)=\EE_x\big[P_{t-s}f(X_s);\ s<\tau_{\cD_k}\big]+\EE_x\big[P_{t-s}f(X_s);\ \tau_{\cD_k}\le s\big].
\]
The first term is continuous in $x$, because the killed process has a continuous transition density. In the second term, the process leaves the ball of radius $\delta$ around $x$ before time $s$; since the drift contributes at most $bs<\delta/2$, the martingale part $Z$ of $X$ must exceed $\delta/2$. By Doob's inequality the second term is at most
\[
 \|f\|_\infty\,\PP\Big(\sup_{u\le s}|Z_u|\ge\tfrac\delta2\Big)\le
 \|f\|_\infty\,\frac{4\,\EE|Z_s|^2}{\delta^2}=\|f\|_\infty\,\frac{8h\operatorname{tr}(\LL)\,s}{\delta^2},
\]
uniformly in $x\in\cK$. Letting $s\downarrow0$, we see that $P_tf$ is the uniform limit of continuous functions on $\cK$, hence continuous.

By Doob's theorem~\cite{DaPratoZabczyk}, \cite[Theorem~1 and Remark~2]{KulikScheutzow}, a strong Feller and irreducible Markov semigroup has at most one invariant probability measure, and $P_t(x,\cdot)$ converges to it in total variation for every $x$. With part (i) this proves the ``if'' direction.

\vspace{0.2cm}

\noindent $\bullet$ Step B ($\Longrightarrow$):~We use the proof by contraposition. Suppose that the condition (C) does not hold. Then, we have the following cases. 

\vspace{0.1cm}

$\diamond$ Step B.1:~Suppose that $G_K$ is disconnected, and let $\mathcal C$ be one of its components. Then $Q=\sum_{a\in\mathcal C}E_a$ is a nonconstant continuous function on $\cD$ and is conserved along every path.

\vspace{0.1cm}

$\diamond$ Step B.2:~Suppose that $G_\Lambda$ is disconnected, and let $\mathcal C$ be one of its components. Then a scalar coordinate of $\sum_{a\in\mathcal C}p_a$ has the same properties and can be used as $Q$. In either case, we choose $q$ so that $\cD_q=\{x\in\cD:Q(x)<q\}$ satisfies $0<\pi_h(\cD_q)<1$, using that $\pi_h$ is strictly positive on $\cD$. Conservation gives $P_t(\mathbf1_{\cD_q}f)=\mathbf1_{\cD_q} P_tf$. Invariance of $\pi_h$ then gives
\[
 \int_{\cD_q} P_tf\,\pi_h\,dx=\int P_t(\mathbf1_{\cD_q}f)\,\pi_h\,dx=\int_{\cD_q} f\,\pi_h\,dx .
\]
Thus the normalised restrictions of $\pi_h$ to $\cD_q$ and to $\cD\setminus\cD_q$ are distinct invariant probability measures. This proves nonuniqueness. Therefore, we have the second assertion.

\vspace{0.1cm}

\noindent (iii)~(Third assertion):~Suppose that $f_0$ and $f_t$ are density functions for the initial law and law of $X_t$ such that 
\[  H(f_0\mid\pi_h)<\infty, \]
and recall that 
\[
\cS_h[f_t]=\int_\cD f_t \,\eta\,dx-h\int_\cD f_t \log f_t \,dx=h\big(\log\cZ_h-H(f_t \mid\pi_h)\big). 
\]
Now, we claim that, for all $t>0$,
\begin{equation} \label{NNN-3}
 H(f_t\mid\pi_h) \leq H(f_0\mid\pi_h)<\infty, \quad \mbox{and the integrals defining $\cS_h[f_t]$ are finite}.
\end{equation}
\noindent $\diamond$ (Monotonicity):~Let $g_0=f_0/\pi_h$, so that $f_t=\pi_h\,P_tg_0$ by Proposition~\ref{f:prop:semigroup}(iii). With $\Phi(u)=u\log u$, which is convex and bounded below, Jensen's inequality for the kernel $P_t(x,\cdot)$ gives 
\[ \Phi(P_tg_0)\le P_t\Phi(g_0), \]
and integrating against the invariant density yields
\[
 H(f_t\mid\pi_h)=\int\Phi(P_tg_0)\,\pi_h\,dx\le\int\Phi(g_0)\,\pi_h\,dx=H(f_0\mid\pi_h)<\infty .
\]
Applied between any two times, this shows that $H(f_t\mid\pi_h)$ is nonincreasing.

\noindent $\diamond$ (Finiteness):~Let $f$ be any probability density on $\cD$ with $H(f\mid\pi_h)<\infty$; by the monotonicity, this applies to $f=f_t$. The function $\eta$ is bounded above on $\cD$, and for $0<\lambda<1/h$ the function $e^{\lambda|\eta|}\pi_h$ is bounded, since $e^{-\lambda\eta}e^{\eta/h}=e^{(1/h-\lambda)\eta}$ where $\eta<0$. The entropy inequality
\[ \int g\,f \,dx\le H(f \mid\pi_h)+\log\int e^{g}\pi_h\,dx \quad \mbox{with $g=\lambda|\eta|$} \]
yields $\int f|\eta|\,dx<\infty$. Since $\log\pi_h=\eta/h-\log\cZ_h$, the integral $\int f\log f\,dx=H(f\mid\pi_h)+\int f\log\pi_h\,dx$ is finite as well. Therefore both integrals in $\cS_h[f_t]$ are finite, the identity $\cS_h[f_t]=h\big(\log\cZ_h-H(f_t\mid\pi_h)\big)$ holds, and $t\mapsto\cS_h[f_t]$ is nondecreasing. This completes the proof.
\begin{remark}
Note that a unique invariant law does not mean that the temperatures along a trajectory become equal: because of the noise they keep fluctuating.
\end{remark}

\subsection{Proof of Theorem \ref{thm:relax}} \label{sec:5.3}
In this subsection, we provide proofs for the assertions in Theorem \ref{thm:relax}. Since the proofs are rather lengthy, we split the proof of each assertion into several subsections.  From now on, we assume {\bf Condition (C)}, so that $R=\LL|_{T\cM}$ is positive definite
\subsubsection{Functional inequalities} \label{sec:5.3.1}
In this part, we provide proofs for the first and second assertions in Theorem \ref{thm:relax}. 

\vspace{0.1cm}

 Note that the rate $\kappa$ in \eqref{eq:kappa} compares the curvature of the entropy with $R$: it is the largest constant such that 
\[ \xi\cdot \mathbb H(x)\xi\ge\kappa\,\xi\cdot R^{-1}\xi \quad \mbox{for all $x\in\cD$ and $\xi\in T\cM$}. \]
We first show that it is positive and give an explicit lower bound.

\begin{lemma}\label{r:lem:hess}
For $x\in\cD$ and variations $\xi_a$ of $p_a$ and $\theta_a$ of $E_a$, one has 
\begin{equation}\label{r:eq:flat}
 -\nabla^2\eta(x)[\xi,\theta]\ \ge\ \frac15\sum_a\Big[\frac{C_a\,\theta_a^2}{E^2}
 +\frac{C_a\,|\xi_a|^2}{\rho_aE}\Big]=:\frac15\,\mathbb H_\flat[\xi,\theta] .
\end{equation}
\end{lemma}
\begin{proof}
We set $A_a=\theta_a-v_a\cdot\xi_a$. Since $T_a\le E/C_a$ and
$|v_a|^2\le 2E/\rho_a$ on $\cD$, formula~\eqref{eq:hessform} gives
\[
 -\nabla^2\eta(x)[\xi,\theta]
 \ge \sum_a\left[\frac{C_a}{E^2}A_a^2+
                    \frac{C_a}{\rho_aE}|\xi_a|^2\right].
\]
On the other hand, $\theta_a=A_a+v_a\cdot\xi_a$, hence
\[
 \theta_a^2\le 2A_a^2+2|v_a|^2|\xi_a|^2
 \le 2A_a^2+\frac{4E}{\rho_a}|\xi_a|^2.
\]
Consequently,
\begin{align*}
 \frac15\left[\frac{C_a\theta_a^2}{E^2}
             +\frac{C_a|\xi_a|^2}{\rho_aE}\right]
 &\le \frac{2C_a}{5E^2}A_a^2
       +\frac{C_a}{\rho_aE}|\xi_a|^2\\
 &\le \frac{C_a}{E^2}A_a^2
       +\frac{C_a}{\rho_aE}|\xi_a|^2.
\end{align*}
Summing over $a$ proves~\eqref{r:eq:flat}.
\end{proof}
\noindent $\bullet$~{\bf (Proof of the first assertion)}:~We claim that the constant $\kappa$ is positive, does not depend on $h$, and
\begin{equation}\label{r:eq:kappabd}
 \kappa\ \ge\ \kappa_\flat:=\frac15\min\Big\{
 \frac{\lambda_2\big(D_{c_v}^{1/2}\cG_\Lambda D_{c_v}^{1/2}\big)}{nE},\
 \frac{\lambda_2\big(D_C^{1/2}\cG_KD_C^{1/2}\big)}{nE^2}\Big\},
\end{equation}
where $D_{c_v}=\operatorname{diag}(c_{v,a})$ and $\lambda_2$ is the second smallest eigenvalue, positive by (C). \newline

\noindent $\diamond$ (Proof of \eqref{r:eq:kappabd}):
Neither $\mathbb H$ nor $R$ depends on $h$. By \eqref{r:eq:flat}, $\kappa$ is at least one fifth of the smallest eigenvalue of $R^{1/2}\mathbb H_\flat R^{1/2}$ on $T\cM$, that is, of $\LL \mathbb H_\flat$ on the invariant subspace $T\cM$. On the whole space, 
\[ \LL \mathbb H_\flat=\frac1{nE}\big(\cG_\Lambda D_{c_v}\otimes I_d\big)\oplus\frac1{nE^2}\cG_KD_C. \]
For a connected weighted graph Laplacian $\cG$ and a positive diagonal matrix $D$, the matrix $\cG D$ is similar to $D^{1/2}\cG D^{1/2}$, its kernel is spanned by $D^{-1}\mathbf 1$, and every eigenvector with a nonzero eigenvalue is orthogonal to $\mathbf 1$, because $\mathbf 1^{\T}\cG D=0$. So the eigenvalues of $\LL \mathbb H_\flat$ on $T\cM$ are the nonzero eigenvalues of the two blocks, and the smallest one is $5\kappa_\flat$. Note that both terms in \eqref{r:eq:kappabd} are rates: since $E=\bar T\sum_bC_b$, the first is of order $\Lambda/(\rho\bar T)$ and the second of order $K/(C\bar T^2)$, as in \eqref{eq:rates}. The bound $\kappa_\flat$ is explicit, but we do not expect it to be sharp. This proves the first assertion. \newline

Next, with the uniform curvature bound we can show the functional inequalities in the second assertion of Theorem \ref{thm:relax}. \newline

\noindent $\bullet$~{\bf (Proof of the second assertion)}:~We claim that for every $h>0$ and every $f\in\cF$,
\begin{equation}\label{r:eq:lsi}
 \Ent_{\pi_h}(f^2)\le\frac2\kappa\,\cE_h(f,f),\quad
 \Var_{\pi_h}(f)\le\frac1\kappa\,\cE_h(f,f).
\end{equation}
\noindent $\diamond$ (Proof of claim \eqref{r:eq:lsi}):~Since $\cF=\overline{C^\infty(\overline\cD)}^{\,\cE_{h,1}}$, we first derive the estimates \eqref{r:eq:lsi} for $f\in C^\infty(\overline\cD)$, and then using density argument, we show that the estimates  \eqref{r:eq:lsi} hold for $f \in \cF$. . \newline 

\noindent $\diamond$~Step A:~Suppose that $f \in C^\infty(\overline\cD)$. Then, we split the verification of  \eqref{r:eq:lsi} in several procedures.
\vspace{0.1cm}

\noindent $\clubsuit$~Step A.1 (Change of variables):~we identify $T\cM$ with $\R^m$ by an isometry, and let $\Omega=\{y\in\R^m:\ x^*+R^{1/2}y\in\cD\}$, which is a bounded convex open set. Then, the image of $\pi_h\,dx$ under $x\mapsto y$ is $\mu(dy)=e^{-U(y)}dy$ with
\[
 U(y)= 
 \begin{cases}
\displaystyle  -\frac1h\,\eta\big(x^*+R^{1/2}y\big)+\text{const}\quad & y\in\Omega, \\
\displaystyle  +\infty \quad & y \notin \Omega.
\end{cases}
\]
For $f\in C^\infty(\overline\cD)$ and $g(y)=f(x^*+R^{1/2}y),$ we have
\[ \nabla_yg=R^{1/2}\nabla f, \quad  |\nabla_yg|^2=\nabla f\cdot\LL\nabla f, \quad \cE_h(f,f)=h\int_\Omega|\nabla g|^2d\mu. \]
On $\Omega$, it follows from the definition of $\kappa$ that
\[ \nabla_y^2U=h^{-1}R^{1/2}\mathbb H R^{1/2}\ge(\kappa/h)\,I. \]

\noindent $\clubsuit$~Step A.2 (Contraction):~Let
\[ s=h/\kappa \quad \mbox{and} \quad \widetilde U(y)=U(y)-|y|^2/(2s). \]
On $\Omega$ the function $\widetilde U$ is smooth with $\nabla^2\widetilde U\ge0$, and $\widetilde U=+\infty$ outside $\Omega$, so $\widetilde U$ is convex on $\R^m$ with values in $(-\infty,+\infty]$. Thus $\mu$ has a density proportional to $e^{-\widetilde U}$ with respect to the Gaussian measure $\gamma_s=N(0,sI)$, and it has compact support. By the Caffarelli contraction theorem~\cite{Caffarelli2000}, in the form of Fathi, Gozlan and Prod'homme~\cite{FathiGozlanProdhomme}, which allows convex functions with the value $+\infty$ and compactly supported targets, applied after the scaling $y\mapsto y/\sqrt s$, there is a $1$-Lipschitz map $S:\R^m\to\R^m$ that pushes $\gamma_s$ forward to $\mu$. Since $S$ is continuous and $\mu$ is supported in $\overline\Omega$, $S$ takes values in $\overline\Omega$. 

\vspace{0.1cm}

\noindent $\clubsuit$~Step A.3 (Transfer):~For Lipschitz $G$, the Gaussian measure $\gamma_s$ satisfies
\[
 \Ent_{\gamma_s}(G^2)\le2s\int|\nabla G|^2d\gamma_s,\quad
 \Var_{\gamma_s}(G)\le s\int|\nabla G|^2d\gamma_s .
\]
For $g$ as in Step A.1, extended smoothly to a neighbourhood of $\overline\Omega$, the function $G=g\circ S$ is Lipschitz and $|\nabla G|\le|\nabla g|\circ S$ almost everywhere. Hence we have
\[
 \Ent_\mu(g^2)=\Ent_{\gamma_s}(G^2)\le2s\int|\nabla g|^2\circ S\,d\gamma_s
 =\frac{2h}\kappa\int_\Omega|\nabla g|^2d\mu=\frac2\kappa\,\cE_h(f,f).
\]
In the same way, we have
\[ \Var_{\pi_h}(f)\le\kappa^{-1}\cE_h(f,f). \]
This verifies \eqref{r:eq:lsi} for $f\in C^\infty(\overline\cD)$. \newline

\noindent $\diamond$~Step B:~Suppose that $f \in \cF$. Then, there exists a sequence $(f_k)$ in $C^\infty(\overline\cD)$ such that
\[ f_k\to f \quad \mbox{in $\cE_{h,1}$}. \]
The variance is continuous in $L^2(\pi_h)$. For the entropy we use the formula 
\[ \Ent_{\pi_h}(g)=\sup \Big \{\int g\,\varphi\,\pi_h\,dx:\ \varphi\ \text{bounded},\quad \int e^{\varphi}\pi_h\,dx\le1 \Big \} \quad \mbox{for $g\ge0$}. \] For each bounded $\varphi$ the map $f\mapsto\int f^2\varphi\,\pi_h\,dx$ is continuous in $L^2(\pi_h)$, so $f\mapsto\Ent_{\pi_h}(f^2)$ is lower semicontinuous, and \eqref{r:eq:lsi} follows. This proves the second assertion. 

\subsubsection{Exponential decay} \label{sec:5.3.2}
In this part, we deal with the exponential relaxation toward an equilibrium. \newline

\noindent $\bullet$~{\bf (Proof of the third assertion)}:~Suppose that 
\[ \eps_a\le1 \quad \mbox{for every $a$}. \]
Then, we claim:
\vspace{0.1cm}
\begin{enumerate}
\item
For $f\in L^2(\pi_h)$, we have
\begin{equation}\label{r:eq:l2}
\Big\|P_tf-\int f\,\pi_h\,dx\Big\|_{L^2(\pi_h)}\le e^{-\kappa t}\,
 \Big\|f-\int f\,\pi_h\,dx\Big\|_{L^2(\pi_h)} .
 \end{equation}
 \item
If the initial law has a density $f_0$ with $H(f_0\mid\pi_h)<\infty$, then
\begin{equation}\label{r:eq:entdecay}
 H(f_t\mid\pi_h)\le e^{-2\kappa t}\,H(f_0\mid\pi_h) .
\end{equation}
\end{enumerate}
Since Condition~(C) holds and $\eps_a\le1$ for every $a$, Theorem~\ref{thm:pos}(i) shows that the solution stays in $\cD$ for all time from every initial state in $\cD$. Hence the standing assumption of Section~\ref{sec:5.1} is satisfied, and Proposition~\ref{f:prop:semigroup} applies.

\noindent {\it Proof of \eqref{r:eq:l2}}:~By Proposition~\ref{f:prop:semigroup}(ii), $P_t$ is the symmetric $L^2(\pi_h)$-semigroup of $(\cE_h,\cF)$. We set 
\[ u_t=P_tf-\int f\,\pi_h\,dx. \]
For $t>0$, $u_t$ is in the domain of the generator and by \eqref{r:eq:lsi},
\[ \frac{d}{dt}\|u_t\|^2=-2\cE_h(u_t,u_t)\le-2\kappa\|u_t\|^2. \]
This yields \eqref{r:eq:l2}. \newline

\noindent {\it Proof of \eqref{r:eq:entdecay}}:~Let $g_0=f_0/\pi_h$. Then,  by Proposition~\ref{f:prop:semigroup}(iii),
\[ f_t=\pi_h\,P_tg_0  \quad \mbox{and} \quad H(f_t\mid\pi_h)=\Ent_{\pi_h}(P_tg_0). \]
Suppose first that $\delta\le g_0\le\delta^{-1}$ for some $\delta>0$. Then $u_t=P_tg_0$ satisfies the same bounds, it is in the domain of $\cL$ for $t>0$, and $\log u_t$ and $\sqrt{u_t}$ belong to $\cF$. With $\Phi(t)=\int u_t\log u_t\,\pi_h\,dx$, Lemma~\ref{f:lem:form}(iv) and the chain rule imply
\begin{align*}
 \Phi'(t)&=\int(\cL u_t)(\log u_t+1)\,\pi_h\,dx=-\cE_h(u_t,\log u_t)\\
 &=-h\int_\cD\frac{\nabla u_t\cdot\LL\nabla u_t}{u_t}\,\pi_h\,dx=-4\,\cE_h\big(\sqrt{u_t},\sqrt{u_t}\big).
\end{align*}
By \eqref{r:eq:lsi} applied to $\sqrt{u_t}$, $\Phi'(t)\le-2\kappa\,\Ent_{\pi_h}(u_t)$. Since $\int u_t\,\pi_h\,dx$ does not depend on $t$, we get 
\[ \frac{d}{dt}\Ent_{\pi_h}(u_t)\le-2\kappa\,\Ent_{\pi_h}(u_t), \]
which gives \eqref{r:eq:entdecay} in this case. \newline

For a general $g_0\ge0$ with $\int g_0\,\pi_h\,dx=1$ and finite entropy, let
\[ g_k=c_k^{-1}\big(\min(g_0,k)+1/k\big), \]
where $c_k$ is the normalizing constant. By dominated convergence, 
\[ g_k\to g_0 \quad \mbox{in $L^1(\pi_h)$} \quad \mbox{and} \quad \Ent_{\pi_h}(g_k) \to \Ent_{\pi_h}(g_0) \quad \mbox{as $k \to \infty$}. \] 
Since $P_t$ is a contraction in $L^1(\pi_h)$, 
\[ P_tg_k\to P_tg_0 \quad \mbox{in $L^1(\pi_h)$}, \]
and the entropy is lower semicontinuous in $L^1(\pi_h)$ by the variational formula used in Step~B in Section~\ref{sec:5.3.1}. Hence, we have
\begin{align*}
 H(f_t\mid\pi_h) &= \Ent_{\pi_h}(P_tg_0)\le\liminf_k\Ent_{\pi_h}(P_tg_k)\\
 &\le e^{-2\kappa t}\Ent_{\pi_h}(g_0) = e^{-2\kappa t}\,H(f_0\mid\pi_h).
\end{align*}
By Theorem~\ref{thm:eq}(iii), the relation \eqref{r:eq:entdecay} also says that the ensemble entropy approaches its maximum $h\log\cZ_h$ at least at the rate $2\kappa$, which does not depend on $h$.

\subsubsection{Linearization and small noise} 
In this part, we provide proofs of the last two assertions in Theorem \ref{thm:relax}. We now compare $\kappa$ with the relaxation rates of the linearized model. At the equilibrium state $x^*$, the gradient $\nabla\eta(x^*)$ is normal to $\cM$, so $\LL\nabla\eta(x^*)=0$. By \eqref{eq:hessform}, $\mathbb H_*:=\mathbb H(x^*)$ is given by
\begin{equation}\label{r:eq:Hstar}
 \xi\cdot \mathbb H_*\xi=\sum_a\Big[\frac{\theta_a^2}{C_a\bar T^2}+\frac{|\xi_a|^2}{\rho_a\bar T}\Big]
 \qquad\text{for }(\xi,\theta)\in T\cM .
\end{equation}

\noindent $\bullet$~{\bf (Proof of the fourth assertion)}:~We claim that 
\vspace{0.2cm}
\begin{enumerate}
\item
The linearisation at $x^*$ of $\dot x=\LL\nabla\eta$ on $\cM$ is 
\[ \dot y=-R\mathbb H_*\,y \quad \mbox{on}~ T\cM. \]
\noindent The operator $R\mathbb H_*$ is diagonalisable, and its eigenvalues are the numbers $\mu^\Lambda_k$, each with multiplicity $d$, and $\mu^K_k$, $k=2,\dots,n$, from \eqref{eq:rates}.
\vspace{0.1cm}
\item
The linearisation of \eqref{eq:rest} at $x^*$ is the Ornstein--Uhlenbeck process:
\[ dy=-R\mathbb H_*\,y\,dt+\sqrt{2hR}\;d\beta \quad \mbox{on $T\cM$}. \] It is reversible with respect to $N(0,h\mathbb H_*^{-1})$. In $L^2$ of this law, the spectrum of minus its generator is the set of sums $\sum_jk_j\mu_j$ with nonnegative integers $k_j$, where $\mu_1,\dots,\mu_m$ are the eigenvalues in (i) counted with multiplicity. Its spectral gap is $\lambda_*$, for every $h>0$.
\end{enumerate}
In the sequel, we provide detailed proofs for the above two assertions. \newline

\noindent (i)~Since $\LL\nabla\eta(x^*)=0$, the linearization is given as follows.
\[ \dot y=\LL\nabla^2\eta(x^*)y=-R\mathbb H_*y. \]
The operator $R\mathbb H_*$ is similar to the symmetric positive definite operator $M:=R^{1/2}\mathbb H_*R^{1/2}$, so it is diagonalisable with positive eigenvalues. By \eqref{r:eq:Hstar}, on the whole space:
\[ \LL(-\nabla^2\eta(x^*))=\frac1{n\bar T}\big(\cG_\Lambda D_\rho^{-1}\otimes I_d\big) \oplus\frac1{n\bar T^2}\cG_KD_C^{-1}, \]
and the argument of \eqref{r:eq:kappabd} shows that its eigenvalues on $T\cM$ are the nonzero eigenvalues of the two blocks, which are the numbers \eqref{eq:rates}.
\vspace{0.2cm}

\noindent (ii)~The noise has constant covariance, so it does not change under linearisation. In the variable $z=R^{-1/2}y$ the equation becomes 
\[ dz=-Mz\,dt+\sqrt{2h}\,d\beta', \]
with a standard Brownian motion $\beta'$ in $T\cM$. In the orthonormal eigenbasis of $M,$ the components are independent one-dimensional Ornstein--Uhlenbeck processes 
\[ dz_j=-\mu_jz_j\,dt+\sqrt{2h}\,d\beta'_j, \]
each reversible with respect to $N(0,h/\mu_j)$. The generator $h\,\partial_j^2-\mu_jz_j\partial_j$ has the Hermite polynomials in $z_j\sqrt{\mu_j/h}$ as eigenfunctions, with eigenvalues $-k\mu_j$, $k=0,1,2,\dots$. Products of these eigenfunctions form an orthogonal basis of $L^2$ of the product law, which is $N(0,hM^{-1})$ in $z$ and $N(0,h\mathbb H_*^{-1})$ in $y$. This gives the spectrum; its smallest nonzero element is $\min_j\mu_j=\lambda_*$. This completes the proof of the second assertion in the above claim.

\vspace{0.5cm}

Note that in the linearized model, $h$ only sets the size of the stationary fluctuations, through the covariance $h\mathbb H_*^{-1}$; it does not change the rates. This is Onsager's regression hypothesis~\cite{Onsager1931} for the linearised model. For the full model we have an upper bound. Now, we are ready to provide the proof of the last assertion in Theorem  \ref{thm:relax}.

\vspace{0.2cm}

\noindent $\bullet$~{\bf (Proof of the fifth assertion)}:~Let $h > 0$ be any. Then, we claim that 
\begin{equation} \label{NNN-6}
 \lambda_{\rm gap}(h)\ge\kappa \quad \mbox{and} \quad \limsup_{h\to0}\lambda_{\rm gap}(h)\le\lambda_*. 
\end{equation}
\noindent $\diamond$ (Derivation of $\eqref{NNN-6}_1$):~By  \eqref{r:eq:lsi}, one has 
\begin{equation} \label{NNN-6-1}
 \kappa \le \frac{\cE_h(f,f)}{\Var_{\pi_h}(f)}. 
\end{equation}
We now take an infimum of the above relation over  $f \in \cF$ with $ \Var_{\pi_h}(f)>0$ to see
\[ \kappa \leq  \lambda_{\rm gap}(h). \]
\noindent $\diamond$ (Derivation of $\eqref{NNN-6}_2$):~Let $\ell\in T\cM$ and $f(x)=\ell\cdot(x-x^*)\in\cF$. Then, we have
\begin{equation} \label{NNN-7}
\cE_h(f,f)=h\,\ell\cdot R\ell. 
\end{equation}
We claim that
\begin{equation}\label{r:eq:laplace}
 \lim_{h\to0}\frac1h\Var_{\pi_h}(f)=\ell\cdot \mathbb H_*^{-1}\ell .
\end{equation}
{\it Proof of \eqref{r:eq:laplace}}:
By Lemma~\ref{r:lem:hess} there is $c_0>0$ with $\mathbb H(x)\ge c_0I$ on $T\cM$ for all $x\in\cD$, and $\nabla\eta(x^*)$ is normal to $\cM$, so we have
\[ \eta(x^*)-\eta(x)\ge\frac{c_0}2|x-x^*|^2 \quad \mbox{on $\cD$}. \]
We set $x=x^*+\sqrt h\,w$. Under $\pi_h$ the variable $w$ has density $\frac{q_h}{\int q_h}$, where
\[
 q_h(w)=\mathbf 1_{\{x^*+\sqrt h\,w\in\cD\}}\,
 \exp\Big(\frac{\eta(x^*+\sqrt h\,w)-\eta(x^*)}h\Big) \quad \mbox{is the unnormalised density}.
\]
It satisfies 
\[ 0\le q_h(w)\le e^{-c_0|w|^2/2} \quad \mbox{and} \quad q_h(w)\to e^{-w\cdot \mathbb H_*w/2} \quad \mbox{for every $w$}, \]
because $x^*$ is an interior point and $\eta$ is smooth there. By dominated convergence, applied to $\int q_h$, $\int w\,q_h$ and $\int w\otimes w\,q_h$, the first and second moments of $w$ converge to those of $N(0,\mathbb H_*^{-1})$. Since $\Var_{\pi_h}(f)=h\,\Var(\ell\cdot w)$, this proves \eqref{r:eq:laplace}. \newline

On the other hand,  we use \eqref{eq:gapdef}, \eqref{NNN-7} and \eqref{r:eq:laplace} to find 
\[
 \limsup_{h\to0}\lambda_{\rm gap}(h)\le\frac{\ell\cdot R\ell}{\ell\cdot \mathbb H_*^{-1}\ell}
 \qquad\text{for every }\ell\in T\cM\setminus\{0\}.
\]
With $\ell=\mathbb H_*^{1/2}u$ the right-hand side is $u\cdot \mathbb H_*^{1/2}R\mathbb H_*^{1/2}u/|u|^2$, whose minimum is the smallest eigenvalue of $\mathbb H_*^{1/2}R\mathbb H_*^{1/2}$. It is similar to $R\mathbb H_*$, so the minimum is $\lambda_*$. This completes the derivation of $\eqref{NNN-6}_2$. 
\begin{remark}
Since $x^*\in\cD$, we also have 
\[ \kappa\le\lambda_{\min}(R^{1/2}\mathbb H_*R^{1/2})=\lambda_*. \]
This inequality can be strict, as the next example shows, so the curvature bound alone cannot give Conjecture~\ref{conj:gap}.
\end{remark}

\subsection{A numerical example}\label{r:sec:num}

To illustrate these bounds in Theorem \ref{thm:relax}, we take $n=2$ and $d=1$. Then $\cM$ is two-dimensional, with coordinates $p=p_1=-p_2$ and $e=E_1$. We take $\rho=(1,1.7)$, $c_v=(1,0.6)$, $\Lambda_{12}=1$, $K_{12}=0.7$ and $E=1$ in arbitrary units. Then $\lambda_*\approx1.604$, a numerical minimisation over $\cD$ gives 
\[ \kappa\approx1.575 \quad \mbox{and} \quad \kappa_\flat\approx0.141. \]
We computed approximations of $\lambda_{\rm gap}(h)$ by the Rayleigh--Ritz method with polynomials in $(p,e)$ of total degree at most $10$, using Gauss--Legendre quadrature on $\cD$. With exact integration, Ritz values are upper bounds for $\lambda_{\rm gap}(h)$; our values are numerical approximations of them, without certified error bounds. Increasing the degree to $12$ does not change the digits in Table~\ref{tab:ritz}.
\begin{table}[ht]
\caption{Ritz approximations of the spectral gap $\lambda_{\rm gap}(h)$ for the example of
Section~\ref{r:sec:num}}\label{tab:ritz}
\centering
\begin{tabular}{@{}c|ccccccc@{}}
\toprule
$h$ & $0.5$ & $0.2$ & $0.1$ & $0.05$ & $0.02$ & $0.01$ & $0.005$\\
\midrule
Ritz value & $2.531$ & $1.983$ & $1.796$ & $1.701$ & $1.643$ & $1.624$ & $1.614$\\
\bottomrule
\end{tabular}
\end{table}
 The values are above $\kappa$ and decrease towards $\lambda_*$ as $h\to0$. This supports Conjecture~\ref{conj:gap}, but it does not prove it.
\section{Galilean covariance}\label{sec:gal}
So far we have worked in the rest frame. In this section we prove Theorem~\ref{thm:gal}, which shows that this is no loss of generality. We first check that the model \eqref{eq:sde} has the Onsager form \eqref{eq:sdeL} in every frame.

\begin{lemma}\label{g:lem:onsager}
The following assertions hold.
\vspace{0.1cm}
\begin{enumerate}
\item
For every state with $T_a>0$ for all $a$,
\begin{equation}\label{g:eq:onsager}
 \LL(v)\nabla\eta=\big((F_a)_a,\ (Q_a+v\cdot F_a)_a\big),
\end{equation}
\item
The noise in \eqref{eq:sde} has covariance $2h\,\LL(v)$.
\end{enumerate}
\end{lemma}
\begin{proof}
\noindent (i)~We write 
\[ \nabla\eta=(\nabla_p\eta,\nabla_E\eta) \quad \mbox{with}~~\nabla_p\eta=(-v_a/T_a)_a \quad \mbox{and} \quad \nabla_E\eta=(1/T_a)_a. \]
Since $A_v^{\T}=\begin{pmatrix}I&B_v^{\T}\\0&I\end{pmatrix}$ and $(B_v^{\T}y)_a=y_a\,v$, we have
\[
 A_v^{\T}\nabla\eta=\big((-v_a/T_a+v/T_a)_a,\ (1/T_a)_a\big)=\big((-u_a/T_a)_a,\ (1/T_a)_a\big).
\]
Applying $\LL=L_{\rm m}\oplus L_{\rm e}$ gives $\big((F_a)_a,(Q_a)_a\big)$, and applying $A_v$ adds $B_vF=(v\cdot F_a)_a$ to the energy component.  \newline

\noindent (ii)~ For the noise, the relations \eqref{eq:noise} gives 
\[ d[M]_t=2hL_{\rm m}\,dt, \quad  d[N]_t=2hL_{\rm e}\,dt \quad \mbox{and} \quad d[M,N]_t=0, \]
so the vector $(dM,\ dN+B_vdM)=A_v(dM,dN)$ has covariance $2hA_v\LL A_v^{\T}=2h\LL(v)$.
\end{proof}

So in a frame with $v\ne0$ the mobility matrix is not block diagonal, and the momentum and energy noises are correlated, with cross covariance $2hL_{\rm m}B_v^{\T}$. The block-diagonal form $\LL$ is correct only in the rest frame.  Now, we are ready to provide detailed proofs for the assertions in Theorem~\ref{thm:gal}. \newline

\subsection{Proof of Theorem~\ref{thm:gal}} 

\noindent $\bullet$~(Proof of first assertion):~The map $G_c$ is affine, and $G_{-c}$ is its inverse. For each $a$,
\[
 \Big(E_a+c\cdot p_a+\tfrac12\rho_a|c|^2\Big)-\frac{|p_a+\rho_ac|^2}{2\rho_a}=E_a-\frac{|p_a|^2}{2\rho_a},
\]
so $G_c$ preserves every $T_a$, and hence $\eta$. The velocities change by $v_a\mapsto v_a+c$ and $v\mapsto v+c$, so $u_a$ does not change. The totals change to $P'=P+\rho c$ and $E'=E+c\cdot P+\tfrac12\rho|c|^2$, and the same computation with $\rho$ and $P$ in place of $\rho_a$ and $p_a$ shows that $\Er$ does not change. Therefore $G_c$ maps $\cD_{P,E}$ onto $\cD_{P',E'}$. The linear part of $G_c$ is $(\xi,\theta)\mapsto(\xi,\theta+B_c\xi)$. It maps $T\cM=\{\sum_a\xi_a=0,\ \sum_a\theta_a=0\}$ to itself, and on $T\cM$ it is unipotent, so its determinant is $1$. Hence $G_c$ also maps Lebesgue measure on $\cM_{P,E}$ to Lebesgue measure on $\cM_{P',E'}$.

\vspace{0.2cm}

\noindent $\bullet$~(Proof of second assertion)):~Summing \eqref{eq:sde} over $a$ shows that $P$ and $E$, and hence $v$, are constant along every path. Let $X$ solve \eqref{eq:sde} and let $X'=G_c(X)$. Since $G_c$ is affine with constant coefficients, It\^o's formula has no second-order term: 
\[ dp_a'=dp_a \quad \mbox{and} \quad dE_a'=dE_a+c\cdot dp_a. \]
By (i), $F_a$ and $Q_a$ are the same functions of $X'$ as of $X$, and $v'=v+c$. Hence, we have
\[
 dp'_a=F_a\,dt+dM_a,\quad
 dE'_a=\big(Q_a+(v+c)\cdot F_a\big)\,dt+dN_a+(v+c)\cdot dM_a ,
\]
which is \eqref{eq:sde} for $X'$ with the same Brownian motions. Taking $c=-v$ gives a solution of the rest-frame model with energy $\Er$.

\vspace{0.2cm}

\noindent $\bullet$~(Proof of third assertion):~The matrix $\LL(v)$ depends on $x$ only through $v=\rho^{-1}\sum_bp_b$, so only derivatives with respect to the momentum coordinates $p_{b,k}$ enter its divergence. In the columns of $p_{b,k}$, the momentum rows contain the constant matrix $L_{\rm m}$, and the energy row of constituent $a$ contains $(B_vL_{\rm m})_{a,(b,k)}=n^{-1}(\cG_\Lambda)_{ab}\,v_k$. Since $\partial v_k/\partial p_{b,k}=1/\rho$,
\[
 (\nabla\cdot\LL)_{E_a}=\sum_{b,k}\frac1{n\rho}(\cG_\Lambda)_{ab}=\frac d{n\rho}\sum_b(\cG_\Lambda)_{ab}=0,
\]
because the rows of a graph Laplacian sum to zero. All other components are zero. For a diffusion with drift $\LL(v)\nabla\eta$ and covariance $2h\LL(v)$, the probability flux of a density $f$ is 
\[ J=\LL(v)\nabla\eta\,f-h(\nabla\cdot\LL)f-h\LL(v)\nabla f. \]
For $f=e^{\eta/h}$ the first and last terms cancel and the middle term is zero, so $J=0$.

\vspace{0.2cm}

\noindent $\bullet$~(Proof of fourth assertion):~Fix $P$ and $E$ with $\Er>0$ and let $c=-v$. By the first and second assertions, $G_c$ maps $\cD_{P,E}$ onto the rest-frame state space $\cD_{0,\Er}$ and maps solutions to solutions. Temperatures are preserved, so Theorem~\ref{thm:pos} holds on $\cD_{P,E}$. Since $G_c$ preserves $\eta$ and Lebesgue measure, the image of the rest-frame Einstein density is again proportional to $e^{\eta/h}$, so Theorem~\ref{thm:eq} holds. The Dirichlet forms, semigroups and functional inequalities are carried over by the same map, so Theorem~\ref{thm:relax} holds with the constants computed in the rest frame with energy $\Er$; in particular $\bar T=\Er/\sum_bC_b$. This completes the proof.

\vspace{0.5cm}

The covariance itself does not need the linear caloric law. If $E_a=\rho_ae_a(T_a)+|p_a|^2/(2\rho_a)$ with $e_a$ smooth and $e_a'>0$, and $\eta=\sum_as_a(T_a)$ with $s_a'(T)=\rho_ae_a'(T)/T$, we still have 
\[ \partial\eta/\partial E_a=1/T_a \quad \mbox{and} \quad \partial\eta/\partial p_a=-v_a/T_a, \]
and Lemma~\ref{g:lem:onsager} and the first three assertions of Theorem~\ref{thm:gal} hold without change. Theorem~\ref{thm:pos}--Theorem \ref{thm:relax}, however, are proved here only for constant specific heats.

\section{Conclusion}\label{sec:disc}
In this paper, we have proposed a stochastic version of the thermomechanical Cucker--Smale model in which the noise acts on the pairwise exchanges of momentum and energy, so that total momentum and total energy are conserved along every path. Requiring the Einstein distribution to be reversible fixes the covariance of the noise and gives each constituent one parameter $\eps_a$, the Boltzmann constant divided by its heat capacity in a cell. For a small noise, this parameter controls the relative temperature variance, with a correction due to the fixed total energy~\cite{LandauLifshitzSP}. When both exchange graphs are connected, all temperatures stay positive for all time if and only if $\eps_a\le1$ for every constituent. For a gas, $\eps_a=2/(D_a\cN_a)$, and this condition is readily satisfied when the cell contains many molecules of every constituent. Under this condition, the Einstein distribution is the unique invariant law, the process converges to it exponentially fast with a convergence rate bounded below independently of the noise strength, and the linearised model has the classical relaxation rates of mixtures.

Several questions still remain open. The threshold $\eps_a=1$ depends on the reference measure (Remark~\ref{rem:ref}), and temperature-dependent exchange coefficients, as suggested by kinetic theory~\cite{BaeHaHwangRuggeri2024}, need a separate analysis near the boundary. Other directions are spatially inhomogeneous interactions, common noise, a derivation of the model from an internal dynamics of the agents, and Conjecture~\ref{conj:gap} for small noise.

\newpage

\begin{appendices}

\section{Proofs of lemmas in Section~\ref{sec:frame}}\label{app:forms}
In this appendix, we provide rigorous proofs for lemmas in Section~\ref{sec:frame} one by one. \newline

\subsection{Proof of Lemma~\ref{f:lem:form}}
\noindent (i) We write $\cE_h(f,f)=h\int_\cD|\LL^{1/2}\nabla f|^2\pi_h\,dx$. Let $f_k$ be smooth functions such that 
\[ f_k \to 0 \quad \mbox{in $L^2(\pi_h)$} \quad \mbox{and} \quad \cE_h(f_k-f_l,f_k-f_l) \to 0. \]
Then we have
\[  \LL^{1/2}\nabla f_k\to w \quad \mbox{in $L^2(\cD,\pi_h)$}. \]
For a smooth vector field $\psi$ with compact support in $\cD$, integration by parts yields
\[ \int_\cD \LL^{1/2}\nabla f_k\cdot\psi\,\pi_h\,dx=-\int_\cD f_k\,\nabla\cdot(\pi_h\LL^{1/2}\psi)\,dx. \]
The divergence is bounded and supported in a compact set where $\pi_h$ is bounded below, so the right side tends to zero. Hence, we have
\[ w=0 \quad \mbox{and} \quad \cE_h(f_k,f_k) \to 0. \]
The boundary $\partial\cD$ is Lebesgue-null, so it carries no mass of $\pi_h\,dx$. For a smooth normal contraction $\varphi$, one has 
\[ \nabla(\varphi\circ f)=\varphi'(f)\nabla f,  \quad  |\varphi'|\le1, \]
so we have
\[ \cE_h(\varphi\circ f,\varphi\circ f)\le\cE_h(f,f), \]
and the closure of a closable Markovian symmetric form is a Dirichlet form~\cite{FOT}. By the Stone--Weierstrass theorem, $C^\infty(\overline\cD)$ is uniformly dense in $C(\overline\cD)$ and $C_c^\infty(\cD)$ in $C_c(\cD)$, so both forms are regular. They are strongly local, since $\cE_h(f,g)=0$ when $f$ is constant near the support of $g$; this holds on the cores and passes to the closures. Finally we have
\[ 1\in\cF \quad \mbox{and} \quad \cE_h(1,1)=0, \]
so $(\cE_h,\cF)$ is recurrent and therefore conservative~\cite{FOT}. \newline

\noindent (ii) Extend $f$ to a Lipschitz function on $\cM$ with the same Lipschitz constant, and let $f_\delta$ be its mollifications. Then we have
\[ f_\delta\to f~\mbox{uniformly,} \quad  \nabla f_\delta\to\nabla f~\mbox{a.e.}, \quad |\nabla f_\delta|~\mbox{is uniformly bounded,} \]
so by dominated convergence theorem, we have
\[ f_\delta\to f \quad \mbox{in $\cE_{h,1}$}. \]  
If $f$ has compact support in $\cD$, we use its zero extension, which is Lipschitz, before mollifying. Then we have
\[ f_\delta\in C_c^\infty(\cD)~~\mbox{for small $\delta$}, \]
and the same argument gives $f\in\cF_0$. \newline

\noindent (iii) For $f\in C_c^\infty(\cD)$ and $g\in C^\infty(\overline\cD)$, the divergence form in \eqref{f:eq:gen} and integration by parts give the identity; there is no boundary term because $f$ has compact support in $\cD$. Both sides are continuous in $g$ for the norm $\cE_{h,1}^{1/2}$, since $\cL f$ is bounded, so the identity extends to $\cF$. \newline

\noindent (iv) Let $\cK\subset\cD$ be compact. On $\cK$, $\pi_h\ge c_\cK>0$, and by {\bf Condition (C)}  we have
\[ \nabla f\cdot\LL\nabla f\ge\lambda^+_{\min}|\nabla f|^2, \]
so 
\[ \int_\cK|\nabla f|^2dx\le(hc_\cK\lambda^+_{\min})^{-1}\cE_h(f,f) \quad \mbox{for smooth $f$}. \]
A Cauchy sequence of smooth functions for $\cE_{h,1}$ is therefore a Cauchy sequence in $H^1$ of the interior of every such $\cK$. Its limit has a weak gradient in $L^2_{\rm loc}(\cD)$, and the formula for $\cE_h$ passes to the limit.

\subsection{Proof of Lemma~\ref{f:lem:ident}}
Let $\cK_j=\{x\in\cD:\operatorname{dist}(x,\partial\cD)\ge1/j\}$. By Lemma~\ref{f:lem:form}(iii) and the Fukushima decomposition~\cite{FOT}, for each $f\in C_c^\infty(\cD)$ the process $f(Y_t)-f(Y_0)-\int_0^t\cL f(Y_s)\,ds$ is a $\PP_y$-martingale for quasi-every $y$, and the same holds for $Y^0$ up to $\zeta$. Using a countable family of such $f$ that determines the martingale problem, and the fact that a countable union of exceptional sets is exceptional, we see that $Y$ and $Y^0$ solve the martingale problem for $\cL$ up to $\sigma_j$ and $\sigma_j^0$, respectively, for quasi-every start. These are the exit times from the open sets $\{x\in\cD:\operatorname{dist}(x,\partial\cD)>1/j\}$, with value zero when the initial point lies outside the corresponding open set. Fix $j$, and replace the coefficients of $\cL$ outside a neighbourhood of $\cK_j$ by bounded, globally Lipschitz ones. The martingale problem for the modified operator is well posed, and a well-posed martingale problem determines the law of any solution up to the exit time from a set where the operators agree~\cite{StroockVaradhan}. So the three stopped processes at level $j$ have the same law. For $i\le j$, the time $\tau_i$ is the first time the path comes within distance $1/i$ of $\partial\cD$, so $\tau_i$ and $X_{\cdot\wedge\tau_i}$ are functions of $X_{\cdot\wedge\tau_j}$, and the same is true for $Y$ and $Y^0$. Hence the joint laws coincide up to every level $j$.

For the second statement, $\{\tau\le t\}=\bigcap_j\{\tau_j\le t\}$ because $\tau_j\uparrow\tau$, and $\{\zeta\le t\}=\bigcap_j\{\sigma^0_j\le t\}$ because $Y^0$ is killed exactly when it leaves every $\cK_j$. For $Y$, if $\sigma_j\le t$ for all $j$, then $\sigma_j$ increases to some $t^*\le t$, and by continuity $Y_{t^*}=\lim_jY_{\sigma_j}$ lies on $\partial\cD$. So $\{\sigma\le t\}=\bigcap_j\{\sigma_j\le t\}$ as well.

\subsection{Proof of Lemma~\ref{f:lem:ac}}
Under (C) the operator $\cL$ is uniformly elliptic on compact subsets of $\cD$, with smooth coefficients. For a smooth domain $\cD'\Subset\cD$ the transition kernel of $X$ killed on leaving $\cD'$ has a continuous density, which is strictly positive on $\cD'\times\cD'$ for $t>0$ when $\cD'$ is connected~\cite{Friedman}. Since $\cD$ is convex, we can choose smooth connected domains $\cD_k\Subset\cD$ increasing to $\cD$. Then $\PP_x(X_s\in A,\ s<\tau)=\lim_k\PP_x(X_s\in A,\ s<\tau_{\cD_k})=0$ for every Lebesgue-null set $A$. $\qed$

\end{appendices}

\backmatter

\bmhead{Acknowledgements}
The work of T. Ruggeri was carried out in the framework of the activities of the Italian National Group for Mathematical Physics (Gruppo Nazionale per la Fisica Matematica, GNFM/INdAM).

\section*{Statements and Declarations}

\noindent\textbf{Funding.}  The work of S.-Y. Ha was supported by National Research Foundation grant funded
by the Korea government (RS-2025-00514472). The work of R. Hu and X. Yang was partially supported by ONR grant N00014-24-1-2432, the Simons Foundation (MP-TSM-00002783), and NSF grant DMS-2420988.

\vspace{0.2cm}

\noindent\textbf{Competing interests.} The authors have no competing interests to declare that are relevant to the content of this article.

\vspace{0.2cm}

\noindent\textbf{Data and code availability.} No datasets were generated or analysed in this study. The scripts for the numerical example in Section~\ref{r:sec:num} are available from the corresponding author on request.

\bibliography{refs}

\end{document}